\documentclass[10pt,epsfig]{amsart}
\usepackage{graphicx,color,epsf}
\input{diagrams}
\newtheorem{theorem}{Theorem}[section]
\newtheorem{lemma}[theorem]{Lemma}
\newtheorem{corollary}[theorem]{Corollary}

\theoremstyle{definition}
\newtheorem{definition}[theorem]{Definition}

\theoremstyle{remark}
\newtheorem{remark}[theorem]{Remark}

\numberwithin{equation}{section}

\newcommand {\hide}[1]{}
\newtheorem{proposition}{Proposition}
 \newcommand {\sign} {\mbox{\bf sign}}

\newcommand {\junk}[1]{}
 \newtheorem{algorithm}{\sc Algorithm}

\newcommand {\R} {\mbox{\rm R}}
\newcommand {\s}        {\mbox{\rm sign}}
\newcommand {\D}     {\mbox{\rm D}}

\newcommand {\C}     {\mbox{\rm C}}

\newcommand {\Real}[1]   {\mbox{${\Bbb R}^{#1}$}}
\newcommand {\Sphere}{\mbox{${\bf S}$}}     
 \newcommand {\re}         {\Real{}}

\newcommand {\Z}  {{\Bbb Z}}

\newcommand{\F}{\mathbb Q}

\newcommand {\RR} {{\mathcal R}}

\newcommand {\Der} {{\rm Der}}

\newcommand {\la}   {{\langle}}
\newcommand {\ra}   {{\rangle}}
\newcommand {\eps} {{\varepsilon}}
\newcommand {\E} {{\rm Ext}}
\newcommand {\dist} {{\rm dist}}

\def\sign{{\rm sign}}

\newcommand {\Ker}      {\mbox{\rm Ker}}

\newcommand {\spanof} {{\rm span}}
\newcommand {\Tot} {{\rm Tot}}

\def\addots{\mathinner{\mkern1mu
\raise1pt\vbox{\kern7pt\hbox{.}}
\mkern2mu\raise4pt\hbox{.}\mkern2mu
\raise7pt\hbox{.}\mkern1mu}}

\begin{document}
\title[Computing Betti Numbers of Sets Defined by Quadratic Inequalities]
{Computing the Top Betti Numbers of Semi-algebraic Sets 
Defined by Quadratic Inequalities in Polynomial Time}
\author{Saugata Basu}
\address{School of Mathematics,
Georgia Institute of Technology, Atlanta, GA 30332, U.S.A.}
\email{saugata@math.gatech.edu}
\thanks{The author was supported in part by an NSF Career Award 0133597 and a 
Sloan Foundation Fellowship.
A preliminary version of this paper appears in the 
Proceedings of the ACM Symposium on Theory of Computing, 2005.}





\begin{abstract}
For any $\ell > 0$, we present an algorithm 
which takes as input a semi-algebraic set, 
$S$, defined by 
$P_1 \leq 0,\ldots,P_s \leq 0$, where each $P_i \in \R[X_1,\ldots,X_k]$
has degree $\leq 2,$
and computes the top $\ell$ Betti numbers of $S$,
$b_{k-1}(S), \ldots, b_{k-\ell}(S),$ in polynomial time. 
The complexity of the algorithm, stated more precisely,  is
$ 
\sum_{i=0}^{\ell+2} {s \choose i} k^{2^{O(\min(\ell,s))}}.
$
For fixed $\ell$, the complexity of the algorithm can be expressed as
$s^{\ell+2} k^{2^{O(\ell)}},$
which is polynomial in the input parameters $s$ and $k$.
To our knowledge this is the first polynomial time algorithm for computing
non-trivial topological invariants of semi-algebraic sets in $\R^k$ 
defined by polynomial inequalities, where
the number of inequalities is not fixed  and the polynomials are
allowed to have degree greater than one.
For fixed $s$, we obtain by letting $\ell = k$, 
an algorithm for computing all
the Betti numbers of $S$ whose complexity is $k^{2^{O(s)}}$.
\end{abstract}

\maketitle
\section{Introduction}
\label{sec:intro}
Let $\R$ be a real closed field and $S \subset \R^k$ a semi-algebraic set
defined by a Boolean formula with atoms of the form
$P > 0, P < 0, P=0$ for $P \in {\mathcal P} \subset \R[X_1,\ldots,X_k]$. 
It is known \cite{O,OP,Milnor,Thom,B99,GV} 
that the topological complexity of $S$ 
(measured by the various Betti numbers of $S$) is bounded by $O(s^2d)^k$,
where $s  = \#({\mathcal P})$ and $d = \max_{P\in {\mathcal P}}{\rm deg}(P).$
Note that these bounds are singly exponential in $k$.
More precise bounds on the individual 
Betti numbers of $S$ appear in \cite{B03}.
Designing efficient algorithms for computing 
the homology groups, and in particular the Betti numbers,  of semi-algebraic
sets are considered amongst the most important problems in algorithmic
semi-algebraic geometry.

Even though the Betti numbers of $S$ are bounded singly exponentially
in $k$, there is no known algorithm for producing a singly exponential
sized triangulation of $S$ (which would immediately imply a singly exponential
algorithm for computing the Betti numbers of $S$). In fact, the
existence of a singly exponential sized triangulation, is considered to
be a major open question in real algebraic geometry. 
Doubly exponential algorithms 
(with complexity $(sd)^{2^{O(k)}}$)
for computing all the Betti numbers are known, 
since it is possible to obtain a triangulation of $S$ in doubly
exponential time using cylindrical algebraic decomposition.
In the absence of singly exponential   algorithms for computing
triangulations of semi-algebraic sets, singly exponential   algorithms
are known only for the problems of testing emptiness \cite{R92,BPR95},
computing the zero-th Betti number 
(i.e. the number of semi-algebraically connected components of  $S$) 
\cite{GV92,HRS94,Canny93a,GR92,BPR99},
as well as the Euler-Poincar\'e characteristic of $S$ \cite{B99}.
Very recently a singly exponential time algorithm has  been given for
computing the first few Betti numbers of semi-algebraic sets \cite{B05} 
(see also \cite{BPR05}).

In this paper, we consider a restricted class of semi-algebraic sets -- namely,
semi-algebraic sets defined by a conjunction of 
quadratic inequalities. Since sets defined
by linear inequalities have no interesting topology, sets  defined
by quadratic inequalities can be considered to be the simplest 
class of semi-algebraic sets which can have non-trivial topology. 
Such sets are in fact quite general, since every semi-algebraic set 
can be
defined by (quantified) formulas involving only quadratic polynomials
(at the cost of increasing the number of variables and the size of the
formula). Moreover, as in the case of general semi-algebraic sets,
the Betti numbers of such sets can be exponentially large.
For example, the set  $S \subset \R^k$ defined by 
\[
X_1(1 - X_1) \leq 0, \ldots, X_k(1 - X_k) \leq 0,
\]
has $b_0(S) = 2^k$.

Hence, it is somewhat surprising that for any fixed constant $\ell$, 
the Betti numbers
$b_{k-1}(S),\ldots,b_{k-\ell}(S)$, of a basic closed semi-algebraic set
$S \subset \R^k$ defined by quadratic inequalities, are  polynomially 
bounded. The following theorem appears in \cite{B03}.

\begin{theorem}
\label{the:quadratic}
Let  $\R$ a real closed field and $S \subset \R^k$ be defined by 
\[
P_1 \leq 0,\ldots, P_s \leq 0, \deg(P_i) \leq 2, 1 \leq i \leq s.
\]
Then, for $\ell \geq 0$,
\[
b_{k-\ell}(S) \leq {s \choose {\ell}} k^{O(\ell)}.
\]
\end{theorem}

Notice that for fixed $\ell$ this gives a polynomial bound on the
highest $\ell$ Betti numbers of $S$ (which could possibly be non-zero). 
Observe also that similar bounds do not hold for sets defined
by polynomials of degree greater than two. For instance, the
set $V \subset \R^k$ defined by the single quartic inequality,
\[
\sum_{i=1}^k X_i^2(X_i - 1)^2 -\varepsilon \geq 0,
\]
will have $b_{k-1}(V) = 2^k$, for all small enough $\varepsilon > 0$.

To see this observe that for all sufficiently small $\varepsilon > 0$,
$\R^k \setminus V $ is defined by 
\[ 
\sum_{i=1}^k X_i^2(X_i - 1)^2 <  \varepsilon.
\]
and has $2^k$ connected components, 
since it retracts onto the set $\{0,1\}^k$. 
It now follows by Alexander duality that 
\[
b_{k-1}(V) = b_0(\R^k \setminus V) = 2^k.
\]

\begin{remark}
Even though Theorem \ref{the:quadratic} is stated for semi-algebraic sets
defined by a conjunction of weak inequalities, 
there is an easy reduction to this case
for basic semi-algebraic sets defined by equalities and strict inequalities.
The same reduction is also applicable to  the algorithmic results described
later in this paper, in particular Theorem \ref{the:main}.

Note that the definition of cohomology groups of 
basic semi-algebraic sets defined by equalities and strict inequalities
over an arbitrary real closed field $\R$ requires some care,
and several possibilities exist. 
In this paper we follow the definition given in \cite{BPR03} which agrees
with singular cohomology in case $\R = {\mathbb R}$. We refer the reader to
\cite{BPR03} for further elaboration on this point.

We have the following easy corollary of Theorem \ref{the:quadratic}.
\begin{corollary}
\label{cor:quadratic}
Let  $\R$ a real closed field and $S \subset \R^k$ be defined by 
\[ 
\bigwedge_{P \in {\mathcal P}_1} P = 0
\bigwedge_{P \in {\mathcal P}_2} P > 0
\bigwedge_{P \in {\mathcal P}_3} P < 0
\]
with $\deg(P) \leq 2$ for each $P \in \bigcup_{i=0,1,2} {\mathcal P}_i$,
and $\# \bigcup_{i=0,1,2} {\mathcal P}_i = s$.

Then, for all $\ell \geq 0$,
\[
b_{k-\ell}(S) \leq {s \choose {\ell}} k^{O(\ell)}.
\]
\end{corollary}

In the following proof, as well as later in the paper,
we will extend the ground field $\R$ by infinitesimal
elements.
We denote by $\R\langle \zeta\rangle$  the real closed field of algebraic
Puiseux series in $\zeta$ with coefficients in $\R$ (see \cite{BPR03} for
more details). 
The sign of a Puiseux series in $\R\langle \zeta\rangle$
agrees with the sign of the coefficient
of the lowest degree term in
$\zeta$. 
This induces a unique order on $\R\langle \zeta\rangle$ which
makes $\zeta$
infinitesimal: $\zeta$ is positive and smaller than
any positive element of $\R$.
When $a \in \R\la \zeta \ra$ is bounded 
from above and below by some elements of $\R$,
$\lim_\zeta(a)$ is the constant term of $a$, obtained by
substituting 0 for $\zeta$ in $a$.
Given a semi-algebraic set
$S$ in ${\R}^k$, the {\em extension}
of $S$ to $\R'$, denoted $\E(S,\R'),$ is
the semi-algebraic subset of ${ \R'}^k$ defined by the same
quantifier free formula that defines $S$.
The set $\E(S,\R')$ is well defined (i.e. it only depends on the set
$S$ and not on the quantifier free formula chosen to describe it).
This is an easy consequence of the transfer principle (see for instance
\cite{BPR03}).

\begin{proof}[Proof of Corollary \ref{cor:quadratic}]
Let $0 < \delta \ll \eps \ll 1$ be infinitesimals. We first replace the set
$S$ by the set $S' \subset \R\la\eps\ra^k$ defined by
$S' = \E(S,\R\la\eps\ra) \cap \bar{B}_k(0,1/\eps)$,
where $\bar{B}_k(0,r)$ denotes the closed ball of radius $r$ centered at
the origin. 
It follows from Hardt's triviality theorem for semi-algebraic mappings
\cite{Hardt}
that $b_i(S) = b_i(S')$ for all
$i \geq 0$. We then replace $S'$ by the set 
$S'' \subset \R\la\eps,\delta\ra^k$ defined by,
\[ 
\bigwedge_{P \in {\mathcal P}_1} P \leq 0 \wedge -P \leq 0
\bigwedge_{P \in {\mathcal P}_2} -P + \delta \leq 0
\bigwedge_{P \in {\mathcal P}_3} -P - \delta \leq 0
\bigwedge \eps^2(X_1^2 + \cdots + X_k^2) - 1 \leq 0.
\]
It follows from Hardt's triviality  again that,
$b_i(S') = b_i(S'')$ for all $i \geq 0$.
Now apply Theorem \ref{the:quadratic}. 
\end{proof}
\end{remark}

\begin{remark}
Another point to note is that the bound in Theorem \ref{the:quadratic} 
depends crucially on the structure of the formula defining the semi-algebraic
set $S$ -- namely that $S$ is defined by a {\em conjunction} of  
polynomial inequalities of the form $P \leq 0$ with ${\rm deg}(P) \leq 2$.
The polynomial bound on the highest Betti numbers no longer holds 
if $S$ is defined by a formula involving quadratic inequalities
but having a different structure -- for instance, a disjunction instead
of a conjunction of such inequalities. 
However, the polynomial bound continues to hold  for
certain small variations of the structure of the formula defining $S$.
For instance, we have the following slight generalization of 
Theorem \ref{the:quadratic} whose proof mimics that 
of Theorem \ref{the:quadratic} and which we omit.

\begin{theorem}
\label{the:quadratic2}
Let  $\R$ a real closed field and $S \subset \R^k$ be defined by 
\[
\bigwedge_{1 \leq i \leq s}
(\bigvee_{1 \leq j \leq m} P_{ij} \leq 0),
\deg(P_{ij}) \leq 2, 1 \leq i \leq s, 1 \leq j \leq m.
\]
Then, for $\ell \geq 0$,
\[
b_{k-\ell}(S) \leq {s \choose {\ell}} k^{O(m\ell)}.
\]
\end{theorem}

For fixed $m$ and $\ell$, the above bound is still polynomial in
$s$ and $k$. 
Similarly, the main result of this paper,
namely Algorithm \ref{alg:general} in Section \ref{sec:general},
can be extended to this situation as well with complexity 
polynomial in $s$ and $k$ for fixed $m$ and $\ell$. We will 
omit the details of this extension.
\end{remark}

Semi-algebraic sets defined by a system of quadratic inequalities 
have a special significance in the theory of computational complexity. 
Even though such sets might seem to be the next simplest class of 
semi-algebraic sets after sets defined by linear inequalities, 
from the point of view of 
computational complexity they represent a quantum leap.
Whereas there exist (weakly) polynomial time algorithms for
solving linear programming, solving quadratic
feasibility problem is provably hard.
For instance, it follows from an easy reduction from the problem of 
testing feasibility of a real quartic equation in many variables, that 
the problem of testing whether a system of quadratic inequalities 
is feasible is  ${\rm NP}_{\rm R}$-complete in the Blum-Shub-Smale
model of computation (see \cite{BCSS}). 
Assuming the input polynomials to have integer
coefficients, the same problem is NP-hard in the classical Turing machine
model, since it is also not difficult to see that the Boolean satisfiability 
problem can be posed as the problem of deciding whether a certain 
semi-algebraic set  defined by  quadratic inequalities is empty or not 
(see Section \ref{sec:hardness}). 
Counting the number of connected components of such sets is even harder.
In fact, we  prove (see Theorem \ref{the:hardness}) that 
for $\ell = O(\log k)$,
computing the $\ell$-th Betti number of 
a basic semi-algebraic set defined by quadratic inequalities in $\R^k$
is $\#$P-hard. 
Note that PSPACE-hardness of the problem of
counting the number of connected components for general semi-algebraic
sets were known before \cite{BC,Reif}, and the proofs of these
results extend easily to the  quadratic case.   
In view of these hardness results,  it is unlikely that there exist polynomial
time algorithms for computing the Betti numbers (or even the first
few Betti numbers) of such a set.
In contrast to these hardness results, the polynomial bound on the top 
Betti numbers of sets defined by quadratic inequalities gives rise to
the possibility that these might in fact be computable in polynomial time. 

In this paper we prove  that for each fixed $\ell > 0$,
the top $\ell$ Betti numbers of basic semi-algebraic 
sets defined by quadratic inequalities are computable in polynomial time. 
We will assume that the polynomials given as input to our algorithms
have coefficients in some ordered domain $\D$ contained in a real
closed field $\R$. We will denote the algebraic closure of $\R$ by $\C$.
By complexity of our algorithms we will mean the number of arithmetic
operations including comparisons in the ring $\D$. When $\D = \Z$,
we will also count the number of bit operations.

The main result of this paper is the following.\\
\noindent {\bf Main Result:} 
We present an algorithm (Algorithm \ref{alg:general} in 
Section \ref{sec:general})
which given a set of $s$ polynomials,
${\mathcal P} = \{P_1,\ldots,P_s\} \subset \R[X_1,\ldots,X_k],$
with ${\rm deg}(P_i) \leq 2, 1 \leq i \leq s,$
computes $b_{k-1}(S), \ldots, b_{k-\ell}(S),$ 
where $S$ is the set defined by $P_1 \leq 0,\ldots,P_s \leq 0$.
The complexity of the algorithm  is
\begin{equation}
\label{eqn:complexity}
\sum_{i=0}^{\ell+2} {s \choose i} k^{2^{O(\min(\ell,s))}}.
\end{equation}
If the coefficients of the polynomials in
${\mathcal P}$ are integers  of bitsizes bounded by
 $\tau$, then the bitsizes of the integers
appearing in the intermediate computations and the output
are bounded by $\tau (sk)^{2^{O(\min(\ell,s))}}.$

To our knowledge this is the first polynomial time algorithm for computing
a non-trivial topological invariant of semi-algebraic sets in $\R^k$ 
defined by polynomial inequalities, where
the number of constraints is not fixed and the polynomials are
allowed to have degree greater than $1$. 
The special case when all the polynomials 
are linear reduces to the well-studied problem of linear programming. 
In this case the set $S$ is either a convex polyhedron or empty, and (weakly) 
polynomial time algorithms are known to decide emptiness of such a set.
In another direction, Barvinok \cite{Barvinok} designed a polynomial
time algorithm for deciding feasibility of systems of quadratic inequalities,
but under the condition that the number of inequalities is bounded by a 
constant (see also \cite{GP} for an interesting generalization as well as a 
constructive version of this result).

\section{Brief Outline}
\label{sec:outline}
Given any compact semi-algebraic set $S$,
we will denote by $b_i(S)$ the rank of $H^i(S,\F)$ (the $i$-th
simplicial cohomology group of $S$ with coefficients in $\F$).
We denote by $\Sphere^k \subset \R^{k+1}$ the unit sphere 
centered at the origin.
We first consider the case of 
semi-algebraic subsets of the unit sphere, $\Sphere^k \subset \R^{k+1},$
defined by homogeneous quadratic inequalities.
We then show how to reduce the general problem to this special case.

Let $S \subset \Sphere^{k}$ be the set defined on $\Sphere^k$  by 
$s$ inequalities, $P_1 \leq 0,\ldots, P_s \leq 0$, where
$P_1,\ldots,P_s \in \R[X_0,\ldots,X_k]$ are homogeneous quadratic
polynomials. For  each $i, 1 \leq i \leq s$, let $S_i \subset \Sphere^k$ denote 
the set defined on $\Sphere^k$ by $P_i \leq 0$. Then, $S = \cap_{i=1}^s S_i.$
There are two main  ingredients in the polynomial time algorithm for
computing the top Betti numbers of $S$.

The first main idea is to consider $S$ as the intersection of the
various $S_i$'s and to utilize  the double complex arising from the generalized
Mayer-Vietoris exact sequence (see Section \ref{sec:top_prelim}). 
It follows from the exactness of the generalized Mayer-Vietoris sequence
(see Proposition \ref{prop:GMV} below),
that the top dimensional homology groups of $S$ are isomorphic to those
of the total complex associated to a suitable truncation of the Mayer-Vietoris
double complex. However, computing even the truncation of the 
Mayer-Vietoris double complex, starting from a triangulation of $S$ would 
entail a doubly exponential complexity. However, we utilize the fact that
terms appearing in the truncated complex depend on
the unions  of the $S_i$'s taken at most $\ell+2$ at a time. 
There are at most $\sum_{j=1}^{\ell+2} {s \choose j}$ such sets.
Moreover, for semi-algebraic
sets defined by the disjunction of a small number of quadratic inequalities, 
we are able to compute in polynomial (in $k$) time a complex,
whose homology groups are  isomorphic to those of the given sets. 
The construction
of these complexes in polynomial time is the second important ingredient
in our algorithm and  is described in detail 
in Section \ref{sec:main} ( Algorithm \ref{alg:union}).
These complexes along with the homomorphisms between them
define
another double complex
whose associated spectral sequence (corresponding
to the column-wise filtration) is isomorphic 
from the $E_2$ term onwards 
to the
corresponding one of the (truncated) Mayer-Vietoris double complex
(see Theorem \ref{the:main} below).
Since, we know that 
the latter converges to the homology groups of $S$, 
the top Betti numbers of $S$ are equal to
the ranks of the homology groups of the associated total complex of the 
double complex we computed. These can then be computed using well known 
efficient algorithms from linear algebra.

The rest of the paper is organized as follows.
In Section \ref{sec:top_prelim} we recall certain basic facts from
algebraic topology including the notions of complexes, and double complexes
of vector spaces, spectral sequences and triangulations of semi-algebraic sets.
We do not prove any results since all of them are quite classical and we
refer the reader to appropriate references \cite{Bredon,Mcleary,BPR03} 
for the proofs.
In Section \ref{sec:algo_prelim}, we recall some basic algorithms in
semi-algebraic geometry that we will need later. We state the inputs, outputs
and complexities of these algorithms.
For classical results in the field of real algebraic geometry,
as well as the details of certain algorithms that we use, we 
refer the reader to \cite{BPR03} for convenience.
In Section \ref{sec:top_quadratic} we describe certain
topological properties of semi-algebraic sets defined by quadratic inequalities
which are crucial for our algorithm. 
Most of  the results in this section are due to Agrachev \cite{Agrachev}.
In Section \ref{sec:comp} we prove the main mathematical results necessary
for our algorithm. 
In Section \ref{sec:main} we describe our  algorithm for computing the top
Betti numbers of semi-algebraic sets defined by quadratic forms. 
We treat the general case in Section \ref{sec:general}.
Finally, in Section \ref{sec:hardness}, we show the computational hardness
of computing the first few Betti numbers 
of a given semi-algebraic set defined by quadratic inequalities,
by proving that the problem is $\#$P-hard.
 
\section{Topological Preliminaries}
\label{sec:top_prelim}
In this section we recall some basic facts from algebraic topology,
related to double complexes, and spectral sequences associated to
double complexes.
We also fix our notations for these objects. All facts that we
need are well known, and we merely give a brief overview, 
referring  the reader to \cite{Bredon, Mcleary} for detailed proofs.

\subsection{Complex of Vector Spaces}
A sequence $\{C^p\}$, $p \in {\mathbb Z}$, of
$\F$-vector spaces
together with a sequence
$\{\partial^p\}$ of
homomorphisms $\partial^p :C^p \rightarrow C^{p+1}$ for
which $\partial^{p+1}\circ \partial^p = 0$ for all
$p$ is called a complex. 

The cohomology groups, $H^p(C^\bullet,\F)$ are defined by,
\[ 
H^p(C^\bullet,\F) = {Z^p(C^\bullet)}/{B^p(C^\bullet)},
\]
where
$B^p(C^\bullet)= {\rm Im}(\partial^{p-1}),$
and 
$Z^p(C^\bullet) = \Ker(\partial^p).$

The cohomology groups, $H^*(C^\bullet,F),$ are all
$\F$-vector spaces
(finite dimensional if the vector spaces $C^p$'s are themselves finite
dimensional). We will henceforth omit reference to 
the field of coefficients $\F$ which is fixed throughout the rest of the
paper.

Given a  complex
$C^\bullet$, we denote by
$\check{C}_\bullet$ the dual complex,
\[
\begin{array}{ccccccc}
\cdots & \longleftarrow & \check{C}_{p-1} &
\stackrel{\check{\partial}_p}{\longleftarrow} & \check{C}_{p} & \longleftarrow &\cdots
\end{array}
\]
where $\check{C}_p = {\rm Hom}(C^p,\F)$ is the vector space dual to $C^p$ and
$$
\check{\partial}_p : {\rm Hom}(C^{p},\F) \rightarrow
{\rm Hom}(C^{p-1},\F)
$$ is the homomorphism dual to $\partial^{p-1}$.

Moreover, 
$H_*(\check{C}_{\bullet},\F) \cong H^*(C^{\bullet},\F).$ 

Given two  complexes, $C^\bullet = (C^p,\partial^p)$ and $D^{\bullet}=
(D^{p},\partial^p)$,
a homomorphism of complexes,
$\phi: C^{\bullet} \rightarrow D^{\bullet},$ is a
sequence of homomorphisms $\phi^p: C^p \rightarrow D^p$ for which
$\partial^p\circ \phi^p = \phi^{p+1}\circ\partial^p$ for all $p.$

In other words, the following diagram
is commutative.
\[
\begin{array}{ccccccc}
\cdots & \longrightarrow & C^p &
\stackrel{\partial^p}{\longrightarrow} & C^{p+1} & \longrightarrow &\cdots
\\ & &
\Big\downarrow\vcenter{\rlap{$\phi^p$}} & &
\Big\downarrow\vcenter{\rlap{$\phi^{p+1}$}} & & \\ \cdots & \longrightarrow
& D^p &
\stackrel{\partial^p}{\longrightarrow} & D^{p+1} & \longrightarrow &\cdots
\end{array}
\]

A homomorphism of complexes,
$\phi: C^{\bullet} \rightarrow D^{\bullet},$ induces homomorphisms,
$\phi^*: H^*(C^{\bullet}) \rightarrow H^*(D^{\bullet}).$
The homomorphism $\phi$ is called a {\em quasi-isomorphism} if the 
homomorphisms $\phi^*$ are isomorphisms.

\subsection{Double Complexes}
\label{double}
\hide{
In this section, we recall  the basic notions of a 
double complex of vector spaces and associated spectral 
sequences.}
A {\em double complex} is a bi-graded vector space,
\[
{C}^{\bullet,\bullet} = \bigoplus_{p,q \in {\mathbb Z}} C^{p,q},
\]
with co-boundary operators
$d : C^{p,q} \rightarrow C^{p,q+1}$ and
$\delta: C^{p,q} \rightarrow C^{p+1,q}$ and such that $d\delta +\delta d = 0$.
We say that ${C}^{\bullet,\bullet}$ is a first quadrant
double complex, if it satisfies
the condition that $C^{p,q} = 0$ if either $p < 0$ or $q < 0$.
Double complexes lying in other quadrants are defined in an analogous manner.
{\small
\begin{diagram}
\vdots &              & \vdots&               &\vdots&              & \\
\uTo^{d} &              & \uTo^{d}&               &\uTo^{d}&              & \\
C^{0,2}  & \rTo^{\delta}& C^{1,2} & \rTo^{\delta} & C^{2,2}& \rTo^{\delta}&\cdots \\
\uTo^{d} &              & \uTo^{d}&               &\uTo^{d}&              & \\
C^{0,1}  & \rTo^{\delta}& C^{1,1} & \rTo^{\delta} & C^{2,1}& \rTo^{\delta}&\cdots \\
\uTo^{d} &              & \uTo^{d}&               &\uTo^{d}&              & \\
C^{0,0}  & \rTo^{\delta}& C^{1,0} & \rTo^{\delta} & C^{2,0}& \rTo^{\delta}&\cdots \\
\end{diagram}
}

We call the complex, ${\rm Tot}^{\bullet}({C}^{\bullet,\bullet}),$
defined by 
\[
{\rm Tot}^n({C}^{\bullet,\bullet}) = \bigoplus_{p+q=n} C^{p,q},
\]
with differential
\[\D^n  = d + \delta: {\rm Tot}^{n}({C}^{\bullet,\bullet}) 
\rightarrow {\rm Tot}^{n+1}({C}^{\bullet,\bullet}),
\]
to be the {\em associated total complex of ${C}^{\bullet,\bullet}$}.
{\small
\begin{diagram}
&&\vdots &              & \vdots&               &\vdots&              & \\
&\luLine&\uTo^{d} &  \luLine            & \uTo^{d}&  \luLine             &\uTo^{d}& \luLine              & \\
\cdots&\rTo^{\delta}& C^{p-1,q+1} &  \rTo^{\delta}&C^{p,q+1}&\rTo^{\delta}&C^{p+1,q+1}&\rTo^{\delta}&\cdots \\
&\luLine&\uTo^{d} &  \luLine            & \uTo^{d}&  \luLine             &\uTo^{d}&        \luLine      & \\
\cdots&\rTo^{\delta}&C^{p-1,q}  & \rTo^{\delta}& C^{p,q} & \rTo^{\delta} & C^{p+1,q}& \rTo^{\delta}&\cdots \\
&\luLine&\uTo^{d} &  \luLine            & \uTo^{d}&   \luLine            &\uTo^{d}&            \luLine  & \\
\cdots&\rTo^{\delta}& C^{p-1,q-1}  & \rTo^{\delta}& C^{p,q-1} & \rTo^{\delta} & C^{p+1,q-1}& \rTo^{\delta}&\cdots \\
&\luLine&\uTo^{d} &   \luLine           & \uTo^{d}& \luLine               &\uTo^{d}&           \luLine   & \\
&&\vdots &              & \vdots&               &\vdots&              & \\
\end{diagram}

}

\subsection{Spectral Sequences}
A {\em spectral sequence} is a sequence of bigraded complexes
$(E_r, d_r: E^{p,q}_r \rightarrow E^{p+r,q-r+1}_r)$ 
(see Figure \ref{fig:spectral}) such that 
the complex $E_{r+1}$ is obtained from $E_r$ by taking its homology
with respect to $d_r$ (that is $E_{r+1} = H_{d_r}(E_r)$).

{\small
\begin{figure}[hbt] 
\begin{center}
\begin{picture}(0,0)%
\includegraphics{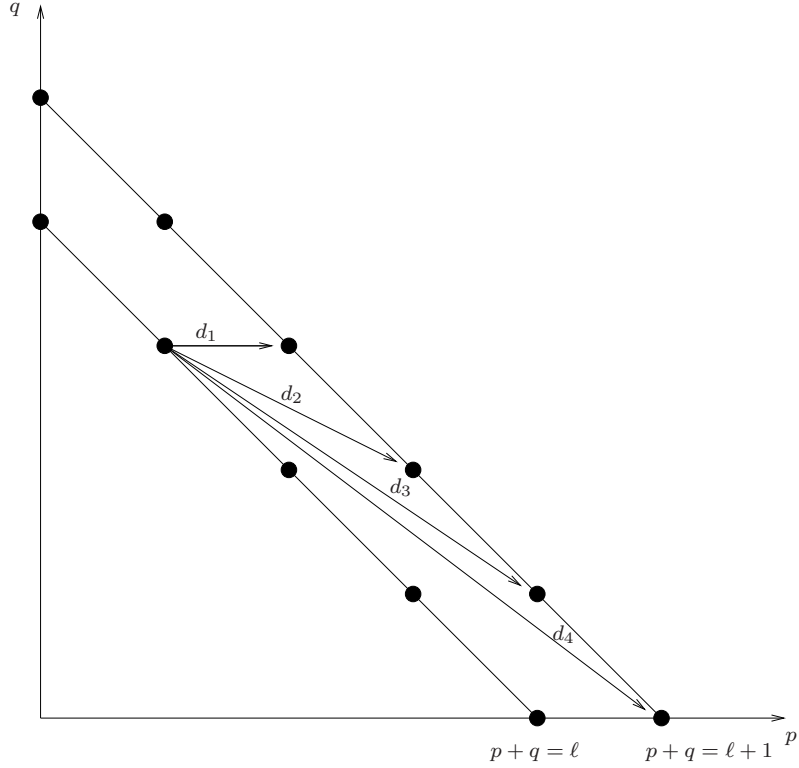}%
\end{picture}%
\setlength{\unitlength}{2565sp}%
\begingroup\makeatletter\ifx\SetFigFont\undefined%
\gdef\SetFigFont#1#2#3#4#5{%
  \reset@font\fontsize{#1}{#2pt}%
  \fontfamily{#3}\fontseries{#4}\fontshape{#5}%
  \selectfont}%
\fi\endgroup%
\begin{picture}(7512,7414)(601,-7394)
\put(6751,-7336){\makebox(0,0)[lb]{\smash{\SetFigFont{8}{9.6}{\familydefault}{\mddefault}{\updefault}{$p + q = \ell+1$}%
}}}
\put(5251,-7336){\makebox(0,0)[lb]{\smash{\SetFigFont{8}{9.6}{\familydefault}{\mddefault}{\updefault}{$p + q = \ell$}%
}}}
\put(8101,-7186){\makebox(0,0)[lb]{\smash{\SetFigFont{8}{9.6}{\familydefault}{\mddefault}{\updefault}{$p$}%
}}}
\put(601,-136){\makebox(0,0)[lb]{\smash{\SetFigFont{8}{9.6}{\familydefault}{\mddefault}{\updefault}{$q$}%
}}}
\put(2401,-3286){\makebox(0,0)[lb]{\smash{\SetFigFont{8}{9.6}{\familydefault}{\mddefault}{\updefault}{$d_1$}%
}}}
\put(3226,-3886){\makebox(0,0)[lb]{\smash{\SetFigFont{8}{9.6}{\familydefault}{\mddefault}{\updefault}{$d_2$}%
}}}
\put(4276,-4786){\makebox(0,0)[lb]{\smash{\SetFigFont{8}{9.6}{\familydefault}{\mddefault}{\updefault}{$d_3$}%
}}}
\put(5851,-6211){\makebox(0,0)[lb]{\smash{\SetFigFont{8}{9.6}{\familydefault}{\mddefault}{\updefault}{$d_4$}%
}}}
\end{picture}
\caption{$d_r: E_r^{p,q} \rightarrow E_r^{p+r, q- r +1}$}
\label{fig:spectral}
\end{center}
\end{figure}
}

There are two spectral sequences,
${'E}_*^{p,q},{''E}_*^{p,q}$,  (corresponding to taking row-wise or
column-wise filtrations respectively) 
associated with a
double complex ${C}^{\bullet,\bullet},$ 
which  will be important for us. 
Both of these converge to $H^*({\rm Tot}^{\bullet}(C^{\bullet,\bullet})).$
This means that the homomorphisms $d_r$ are eventually zero, and hence the
spectral sequences stabilize, and
\[
\bigoplus_{p+q = i} {'E}_{\infty}^{p,q} \cong 
\bigoplus_{p+q = i} {''E}_{\infty}^{p,q} \cong 
H^i({\rm Tot}^{\bullet}(C^{\bullet,\bullet})),
\]
for each $i \geq 0$.

The first terms of these are
${'E}_1 = H_{\delta}(C^{\bullet,\bullet}), 
{'E}_2 = H_dH_{\delta}(C^{\bullet,\bullet})$, and
${''E}_1 = H_d (C^{\bullet,\bullet}), 
{''E}_2 = H_\delta H_d (C^{\bullet,\bullet})$.

Given  two (first quadrant) double complexes, $C^{\bullet,\bullet}$ and
$\bar{C}^{\bullet,\bullet},$ a homomorphism of double complexes,
\[
\phi: C^{\bullet,\bullet} \longrightarrow \bar{C}^{\bullet,\bullet},
\]
is a collection of homomorphisms,
$\phi^{p,q}: C^{p,q} \longrightarrow \bar{C}^{p,q},$ such that
the following diagrams commute.

\[
\begin{array}{ccc}
C^{p,q} &
\stackrel{\delta}{\longrightarrow} & C^{p+1,q} \\ 
\Big\downarrow\vcenter{\rlap{$\phi^{p,q}$}} & &
\Big\downarrow\vcenter{\rlap{$\phi^{p+1,q}$}} \\ 
\bar{C}^{p,q} &
\stackrel{\delta}{\longrightarrow} & \bar{C}^{p+1,q}
\end{array}
\]

\[
\begin{array}{ccc}
 C^{p,q} &
\stackrel{d}{\longrightarrow} & C^{p,q+1}\\
\Big\downarrow\vcenter{\rlap{$\phi^{p,q}$}} & &
\Big\downarrow\vcenter{\rlap{$\phi^{p,q+1}$}} \\ 
\bar{C}^{p,q} &
\stackrel{d}{\longrightarrow} & \bar{C}^{p,q+1}
\end{array}
\]

A homomorphism of double complexes,
\[
\phi: C^{\bullet,\bullet} \longrightarrow \bar{C}^{\bullet,\bullet},
\]
induces homomorphisms,
$\phi_s: E_s \longrightarrow \bar{E}_s$
between the associated spectral sequences (corresponding either to the
row-wise or column-wise filtrations).
For the precise definition of homomorphisms of spectral sequences,
see \cite{Mcleary}.
We will use the following useful fact 
(see \cite{Mcleary}, page 66, Theorem 3.4) several times in the 
rest of the paper.
\begin{theorem}
\label{the:spectral}
If $\phi_s$ is an isomorphism for some $s \geq 1$, then
$E_r^{p,q}$ and $\bar{E}_r^{p,q}$ are isomorphic for
all $r \geq s$. In other words, the induced homomorphism,
$\phi: \Tot^{\bullet}(C^{\bullet,\bullet}) \rightarrow 
\Tot^{\bullet}(\bar{C}^{\bullet,\bullet})$ is a quasi-isomorphism.
\end{theorem}

\subsection{Triangulation of semi-algebraic sets}
A triangulation
of a compact semi-algebraic set $S$ is a simplicial complex 
$\Delta$ together with a
semi-algebraic homeomorphism from $\vert \Delta\vert$ to $S$. 
Given such a triangulation we will often identify the simplices in
$\Delta$ with their images in $S$ under the given homeomorphism, and
will  refer to the triangulation by $\Delta$.

Given a triangulation $\Delta$, the cohomology groups 
$H^i(S)$ are isomorphic to the simplicial cohomology groups 
$H^i(\Delta)$ of the simplicial complex  $\Delta$ and
are in fact independent of the triangulation $\Delta$ 
(this fact is classical over $\re$;  see for instance \cite{BPR03} 
for a self-contained proof in the category of semi-algebraic sets).
 
We call a triangulation $h_1: |\Delta_1| \rightarrow S$
of a semi-algebraic set $S$, to be a {\em refinement}
of a triangulation
$h_2: |\Delta_2| \rightarrow S$ if
for every simplex $\sigma_1 \in \Delta_1$, there exists a simplex
$\sigma_2 \in \Delta_2$ such that $h_1(\sigma_1) \subset h_2(\sigma_2).$

Let $S_1 \subset S_2$ be two compact semi-algebraic subsets of
$\R^k$. We say that a semi-algebraic 
triangulation $h: |\Delta| \rightarrow S_2$ of $S_2$, { \em respects}
$S_1$ if for every simplex $\sigma  \in \Delta$,
$h(\sigma) \cap S_1 = h(\sigma)$ or $\emptyset$.
In this case, $h^{-1}(S_1)$ is identified with a sub-complex
of $\Delta$ and
$h|_{h^{-1}(S_1)} :h^{-1}(S_1) \rightarrow S_1$ is a semi-algebraic 
triangulation  of $S_1$. We will refer to this sub-complex  by 
$\Delta|_{S_1}$.

We will need the following theorem which can be deduced 
from Section 9.2 in  \cite{BCR} (see also \cite{BPR03}).

\begin{theorem}
\label{the:triangulation}
Let $S_1 \subset S_2 \subset\R^k$ 
be closed and bounded semi-algebraic sets, and
let $h_i: \Delta_i \rightarrow S_i, i = 1,2$ be 
semi-algebraic triangulations of
$S_1,S_2$. Then, there exists a semi-algebraic triangulation
$h: \Delta \rightarrow S_2$ of $S_2$, such that
$\Delta$ respects $S_1$,
$\Delta$ is a refinement of $\Delta_2$, and
$\Delta|_{S_1}$ is a refinement of $\Delta_1$.

Moreover, there exists an
algorithm which computes such a triangulation 
whose complexity is bounded by
$(sd)^{2^{O(k)}}$, where $s$ is the number of polynomials used in the
definition of $S_1$ and $S_2$, and $d$ a bound on their degrees. 
\end{theorem}

\subsection{Mayer-Vietoris Double Complex}
\label{subsec:MV}
Let $S_1,\ldots,S_s \subset \R^k$ be closed and bounded semi-algebraic sets,
and let $S = \cap_{i=1}^{s} S_i$.
Choose a sufficiently fine triangulation of $\cup_{1=1}^{s} S_i$,
such that all intersections of the form, 
$S_{i_0} \cap \cdots \cap S_{i_\ell}$, correspond to subcomplexes
of the simplicial complex of the triangulation. Note that the existence
of such a triangulation (in fact, a semi-algebraic triangulation) is well
known (see \cite{BPR03}). However, the best algorithm for computing
such a triangulation has complexity which is doubly exponential 
(in $k$) and produces
a doubly exponential sized triangulation, and hence is not suitable 
for our purpose.

The first main ingredient for our polynomial time algorithm is the double
complex associated to  the generalized Mayer-Vietoris sequence, which
we describe first.

For each $S_i$, let $A_i$ be the subcomplex corresponding to it,
and let $A = A_1 \cap \cdots \cap A_s$. For any $p \geq 0$, let 
$A^{\alpha_0,\ldots,\alpha_p}$ denote the union,
$A_{\alpha_0} \cup \cdots \cup A_{\alpha_p}$.
Let $C_i(A)$ denote the
${\F }$-vector space of $i$-chains of $A$, and 
$C_{\bullet}(A)$ the corresponding chain complex.

\begin{proposition}
\label{prop:GMV}
The following sequence is exact.
$$
\displaylines{
0 \longrightarrow C_{\bullet}(A) \stackrel{i}{\longrightarrow}
\bigoplus_{\alpha_0} C_{\bullet}(A^{\alpha_0}) 
\stackrel{\delta}{\longrightarrow}
\bigoplus_{\alpha_0 < \alpha_1} C_{\bullet}(A^{\alpha_0,\alpha_1}) 
\stackrel{\delta}{\longrightarrow}\cdots\cr
\stackrel{\delta}{\longrightarrow}
\bigoplus_{\alpha_0 < \cdots < \alpha_p}C_{\bullet}(A^{\alpha_0,\ldots,\alpha_p})
\stackrel{\delta}{\longrightarrow}
\bigoplus_{\alpha_0 < \cdots < \alpha_{p+1}}C_{\bullet}(A^{\alpha_0,\ldots,\alpha_{p+1}})\stackrel{\delta}{\longrightarrow}\cdots,
}
$$
where $i$ is induced by inclusion
and the connecting homomorphisms $\delta$ are defined as follows: \\
for $c \in \oplus_{\alpha_0 < \cdots < \alpha_p}
C_{\bullet}(A^{\alpha_0,\ldots,\alpha_p})$, 
$
(\delta c )_{\alpha_0,\ldots,\alpha_{p+1}}
= \sum_{0 \leq i \leq p+1} (-1)^i 
c_{\alpha_0,\ldots,\hat{\alpha_i},\ldots,\alpha_{p+1}}.
$
\end{proposition}
\begin{proof}
See \cite{B03}.
\end{proof}

We have a corresponding (fourth quadrant)
double complex ${\mathcal N}^{\bullet,\bullet},$ 
with
\[
{\mathcal N}^{p,q} 
=  \oplus_{\alpha_0,\ldots,\alpha_p}C_{q}(A^{\alpha_0,\ldots,\alpha_p}).
\]

{\small
\[
\begin{array}{ccccccc}
& & 0 && 0 && 0  \\

& &
\Big\downarrow\vcenter{\rlap{$\partial$}} & &
\Big\downarrow\vcenter{\rlap{$\partial$}} & &
\Big\downarrow\vcenter{\rlap{$\partial$}}   \\

0 & \longrightarrow & \oplus_{\alpha_0}C_{k}(A^{\alpha_0}) &
\stackrel{\delta}{\longrightarrow} & \oplus_{\alpha_0<\alpha_1}
C_{k}(A^{\alpha_0,\alpha_1}) &
\stackrel{\delta}{\longrightarrow} & \oplus_{\alpha_0<\alpha_1 <\alpha_2}
C_{k}(A^{\alpha_0,\alpha_1,\alpha_2}) 
\\

& &
\Big\downarrow\vcenter{\rlap{$\partial$}} & &
\Big\downarrow\vcenter{\rlap{$\partial$}} & &
\Big\downarrow\vcenter{\rlap{$\partial$}}   \\

0 & \longrightarrow & \oplus_{\alpha_0}C_{k-1}(A^{\alpha_0}) &
\stackrel{\delta}{\longrightarrow} & \oplus_{\alpha_0<\alpha_1}
C_{k-1}(A^{\alpha_0,\alpha_1}) &
\stackrel{\delta}{\longrightarrow} & \oplus_{\alpha_0<\alpha_1 <\alpha_2}
C_{k-1}(A^{\alpha_0,\alpha_1,\alpha_2}) 
\\

& &
\Big\downarrow\vcenter{\rlap{$\partial$}} & &
\Big\downarrow\vcenter{\rlap{$\partial$}} & &
\Big\downarrow\vcenter{\rlap{$\partial$}}   \\

0 & \longrightarrow & \oplus_{\alpha_0}C_{k-2}(A^{\alpha_0}) &
\stackrel{\delta}{\longrightarrow} & \oplus_{\alpha_0<\alpha_1}
C_{k-2}(A^{\alpha_0,\alpha_1}) &
\stackrel{\delta}{\longrightarrow} & \oplus_{\alpha_0<\alpha_1 <\alpha_2}
C_{k-2}(A^{\alpha_0,\alpha_1,\alpha_2})
\\

& &
\Big\downarrow\vcenter{\rlap{$\partial$}} & &
\Big\downarrow\vcenter{\rlap{$\partial$}} & &
\Big\downarrow\vcenter{\rlap{$\partial$}}   \\

0 & \longrightarrow & \oplus_{\alpha_0}C_{k-3}(A^{\alpha_0}) &
\stackrel{\delta}{\longrightarrow} & \oplus_{\alpha_0<\alpha_1}
C_{k-3}(A^{\alpha_0,\alpha_1}) &
\stackrel{\delta}{\longrightarrow} & \oplus_{\alpha_0<\alpha_1 <\alpha_2}
C_{k-3}(A^{\alpha_0,\alpha_1,\alpha_2})
\\

& &
\Big\downarrow\vcenter{\rlap{$\partial$}} & &
\Big\downarrow\vcenter{\rlap{$\partial$}} & &
\Big\downarrow\vcenter{\rlap{$\partial$}}   \\

& & \vdots  && \vdots  && \vdots   \\
\end{array}
\]
}

Standard facts about the Mayer-Vietoris spectral sequence then yields
the following well known proposition (see \cite{B03}).
\begin{proposition}
\label{prop:MV}
For $0 \leq i \leq k,$
\[
H_i(A) \cong H^i({\rm Tot}^{\bullet}({\mathcal N}^{\bullet,\bullet})).
\]
\end{proposition}
Note that since ${\mathcal N}^{\bullet,\bullet}$ is a fourth-quadrant
double complex, the complex
${\rm Tot}^{\bullet}({\mathcal N}^{\bullet,\bullet})$ is defined by,
\[
{\rm Tot}^{i}({\mathcal N}^{\bullet,\bullet}) = 
\bigoplus_{p+ k -q = i} {\mathcal N}^{p,q}.
\]

Moreover, if we  denote by ${\mathcal N}_{\ell}^{\bullet,\bullet}$ 
the truncated complex
defined by,
$$
\begin{array}{cccc}
{\mathcal N}_{\ell}^{p,q} &=& {\mathcal N}^{p,q}, &  \;\;
0 \leq p \leq \ell+1,\;\;  k-\ell-1 \leq q \leq k,  \\ 
                    &=& 0,             & \;\;\mbox{otherwise},
\end{array}
$$
then it is clear that, 
\begin{equation}
\label{eqn:tot}
H_i(A) \cong H^i({\rm Tot}^{\bullet}({\mathcal N}_\ell^{\bullet,\bullet})), 
\;\;k-\ell \leq i \leq k. 
\end{equation}

As noted previously, we cannot hope to compute even the truncated complex
${\mathcal N}_{\ell}^{\bullet,\bullet}$ since we do not know how to compute
triangulations efficiently. We overcome this problem by 
computing  another double complex
${\mathcal D}_{\ell}^{\bullet,\bullet}$,
such that there exists a homomorphism of double complexes,
\[
\psi: {\mathcal D}_{\ell}^{\bullet,\bullet} \rightarrow 
{\mathcal N}_{\ell}^{\bullet,\bullet},
\]
which induces an isomorphism between the $E_1$ terms of the
spectral sequences associated to 
the double complexes ${\mathcal D}_{\ell}^{\bullet,\bullet}$ and 
${\mathcal N}_{\ell}^{\bullet,\bullet}.$ 
This implies, by virtue of Theorem \ref{the:spectral}, 
that the homology groups of the associated total complexes are isomorphic,
that is,
\[
H^*({\rm Tot}^{\bullet}({\mathcal N}_{\ell}^{\bullet,\bullet}))
\cong 
H^*({\rm Tot}^{\bullet}({\mathcal D}_{\ell}^{\bullet,\bullet})).
\]

The construction of the double complex 
${\mathcal D}_{\ell}^{\bullet,\bullet}$ is  
described in the Section \ref{sec:comp}.

\section{Algorithmic Preliminaries}
\label{sec:algo_prelim}
In this section, we recall some basic algorithms from semi-algebraic
geometry that we will need later in the paper. Our main reference is
\cite{BPR03}. Here we just describe the inputs, outputs and the complexities
of these algorithms.

\subsection{Sign Conditions}
Let $\R$ be a real closed field and $\C$ its algebraic closure.
 
A {sign condition} is an element of $\{0,1,- 1\}$.
We define
$$
\begin{cases}
{\rm sign}(x) =0& \mbox{if and only if } x=0,\cr
{\rm sign}(x) =1& \mbox{if and only if } x > 0,\cr
{\rm sign}(x) =-1& \mbox{if and only if } x < 0.
\end{cases}
$$

Let  ${\mathcal Q}\subset\R[X_1,\ldots,X_k]$.
A  {sign condition}  on
$\mathcal Q$ is an element of $\{0,1,- 1\}^{\mathcal Q}$.
We say that $\mathcal Q$ {realizes} the sign condition
$\sigma$ at $x\in\R^k$ if
$$\bigwedge_{Q \in {\mathcal Q}}\s(Q(x))=\sigma(Q).$$
The {realization of the sign condition
$\sigma$}  is
$$\RR(\sigma)= \{x\in \R^k\;\mid\;\bigwedge_{Q\in{\mathcal Q}}
                        \s({Q}(x))=\sigma(Q) \}.$$
The sign condition $\sigma$
is {realizable}
if $\RR(\sigma)$ is non-empty.
We denote by $\sign({\mathcal Q})$ the set of realizable sign conditions
of ${\mathcal Q}.$

\subsection{Representations of points}
Since our algorithms will have to deal with points whose
co-ordinates are 
not rational numbers,
we first describe the particular 
representations of such points that we will use.
In the following $\D$ denotes an ordered domain contained in the 
real closed field $\R$.

Let $P \in \R[X]$ and $\sigma \in \{0,1,-1\}^{\Der(P)}$,
a sign condition on the  set $\Der(P)$ of derivatives of $P$.
The sign condition $\sigma$
is  a Thom encoding  of ${x \in \R}$
if $\sigma(P)=0$ and   $\sigma$  is the sign condition taken by
the set $\Der(P)$  at $x$.
A real root $x$ of $P$ is completely 
characterized by the signs of the derivatives
of $P$ at $x$ and
we say that $x$  is specified by $\sigma.$
Given a Thom encoding $\sigma$, we  denote by  $x(\sigma)$
the root of $P$ in $\R$ specified by $\sigma$.

The representations of points are as follows.
A {$k$-univariate representation}
is a $k+2$-tuple of polynomials of
$\D[T]$,
$$(f(T),g_0(T),g_1(T),\ldots ,g_k(T)),
$$
such that $f$ and $g_0$ are coprime.
Note that $g_0(t)\not=0$
if $t\in\C$ is a root of $f(T)$.
The {points associated}
to a univariate
representation are the points $$
\Bigg( {g_1(t) \over g_0(t)},\ldots, {g_k(t)
\over g_0(t)}\Bigg)\in \C^k $$ where
$t\in\C$ is a root of $f(T)$.

A real $k$-univariate representation
a pair $u,\sigma$ where
$u$ is a $k$-univariate representation and $\sigma$ is the Thom
 encoding of a root of $f$,
$t_\sigma \in \R.$
The {point associated} to the real univariate representation is the point
$$
\Bigg( {g_1(t_\sigma) \over g_0(t_\sigma)},\ldots, {g_k(t_\sigma) \over
g_0(t_\sigma)}\Bigg)\in
\R^k.
$$

\hide{
\subsection{Computing Sample Points}
The following algorithm computes a finite set of sample points
(represented by real univariate representations),
which is guaranteed to meet every connected component of the realizations 
of every
realizable sign conditions of a family of polynomials. Its complexity
is polynomial in the number and degrees of the input polynomials and
exponential in the number of variables.

\begin{algorithm} [Sampling] 
\label{13:alg:sampling2}
\item[]
\begin{description}
\item [{\sc Input}] a set of $s$ polynomials,
$$
{\mathcal P} = \{ P_1,\ldots,P_s \} \subset \D[X_1,\ldots,X_k],
$$
each of degree at most $d$.
 \item  [{\sc Output}] a set ${\mathcal U}$ of
real univariate representations in $\D[T]^{k+2}$ such that
 the
associated points  form a set of sample points
for ${\mathcal P}$ in $\R^k$,
meeting every semi-algebraically connected
 component of
$\RR(\sigma)$
for every $\sigma \in {\rm Sign}({\mathcal P})$,
and the  signs of  the elements of ${\mathcal P}$
at these points.
\end{description}
\end{algorithm}

\vspace{.1in}
\noindent
{\sc Complexity:}
The complexity of the algorithm
$$\displaylines{
s\sum_{j \le k}4^j {s  \choose j} d^{O(k)}=
s^{k+1} d^{O(k)}.}
$$ 
However, the number of
points actually constructed is only
$$\displaylines
{\sum_{j \le k}4^j {s  \choose j} O(d)^k.}
$$
}

\subsection{Complexity of Linear Algebra Subroutines}
In our algorithms, we will also have to compute eigen-vectors of,
as well as perform Gram-Schmidt orthogonalizations on,
real symmetric matrices whose entries are 
described by univariate representations.
Suppose that $M$ is a matrix of size $k \times k$, given as
a  real $k^2$-univariate representation.
Recall that this means that $M$ is given as a pair,
$\sigma,u$, where $\sigma$ is a Thom encoding and
$u= (f(T),g_0(T), g_{1,1}(T),\ldots,g_{k,k}(T))$ is a 
$k^2$-univariate representations.

If the degrees of the polynomials $f,g_0,g_{1,1},\ldots,g_{k,k}$ are 
bounded by $D$, 
then using the most naive algorithms, all the above tasks can be performed 
with complexity bounded by $(kD)^{O(1)}.$ 
In order to see this notice that the complexity of the algorithms
of linear algebra with matrices of size $k \times k$ is bounded by
$k^{O(1)}$. However,
since the computations takes
place in the ring $\D[T]/(f)$, each arithmetic operation in this
ring costs $D^{O(1)}$ operations in the ring $\D$.

\subsection{Complexity of Computing  Linear Arrangements}
\label{subsec:arrangement}
We will also need to compute descriptions of cell
complexes, whose cells are the chambers of different dimensions 
in an arrangement of $\ell$ hyperplanes in $\R^k$. Notice that the
number of cells in such an arrangement can be bounded naively by
$3^{\ell}$. 
Using standard algorithms, descriptions of the  cells in such a complex,
and their incidence relationships, can be computed with complexity 
$k^{O(\ell)}$. 

Moreover, in case the arrangement is parametrized, that is the
hyperplanes are described by linear equations, whose coeffcients are 
described by parametrized univariate representations,
with at most $\ell$ parameters, and whose degrees are bounded by $D$,
one can compute a family of polynomials in the parameters, such that
for all values of the parameters 
satisfying a fixed sign condition on this family,
the combinatorial structure of the cell complex remain constant.
Moreover, the complexity of this procedure, as well as the number
and degrees of this new family of polynomials, are all bounded by
$(kD)^{O(\ell)}$.

\section{Topology of sets defined by quadratic constraints}
\label{sec:top_quadratic}
The results of this section were proved by Agrachev \cite{Agrachev} in the 
context of {\em non-degenerate} (see \cite{Agrachev} for the precise 
definition of non-degeneracy) quadratic maps. However, for the
purposes of this paper the assumption of non-degeneracy is not required
as shown below.

Let $P_1,\ldots,P_{s}$ be homogeneous quadratic polynomials in
$\R[X_0,\ldots,X_k]$. 

We denote by $P = (P_1,\ldots,P_s): \R^{k+1} \rightarrow \R^s$,
the map defined by the polynomials $P_1,\ldots,P_s.$  

Let
$$
T = \bigcup_{1 \leq i \leq s}\{ x \in \Sphere^k \mid   P_i(x) \leq   0 \}.
$$ 

Let 
$$
\Omega = \{\omega \in \R^{s} \mid  |\omega| = 1, \omega_i \leq 0, 1 \leq i \leq s\}.
$$

For $\omega \in \Omega$ we denote by ${\omega}P$ the quadratic form
defined by 
\[
{\omega}P = \sum_{i=1}^{s} \omega_i P_i.
\]

Let $B \subset \Omega \times \Sphere^k$ be the set defined by,
\[
B = \{ (\omega,x)\mid \omega \in \Omega, x \in \Sphere^k \;\mbox{and} \; 
{\omega}P(x) \geq 0\}.
\]

We denote by $\phi_1: B \rightarrow \Omega $ and 
$\phi_2: B \rightarrow \Sphere^k$ the two projection maps. 

\begin{diagram}
&&B&& \\
& \ldTo^{\phi_1} &&\rdTo^{\phi_2} & \\
\Omega &&&& \Sphere^k
\end{diagram}

With the notation developed above,
\begin{proposition}
\label{prop:homotopy2}
The map $\phi_2$ gives a homotopy equivalence between $B$ and 
$\phi_2(B) = T$.
\end{proposition}

\begin{proof}
We first prove that $\phi_2(B) = T.$
If $x \in T,$ 
then there exists some $i, 1 \leq i \leq s,$ such that
$P_i(x) \leq 0.$ Then for $\omega = (-\delta_{1,i},\ldots,-\delta_{s,i})$
(where $\delta_{ij} = 1$ if $i=j$, and $0$ otherwise),
we see that $(\omega,x) \in B.$
Conversely,
if $x \in \phi_2(B),$ then there exists 
$\omega = (\omega_1,\ldots,\omega_s) \in \Omega$ such that,
$\sum_{i=1}^s \omega_i P_i(x) \geq 0$. Since, 
$\omega_i \leq 0, 1\leq i \leq s,$ and not all $\omega_i = 0$,
this implies that $P_i(x) \leq 0$ for
some $i, 1 \leq i \leq s$. This shows that $x \in T$.

For $x \in \phi_2(B)$, the fibre 
$$
\phi_2^{-1}(x) = \{ (\omega,x) \mid  
 \omega \in \Omega \;\mbox{such that} \;  {\omega}P(x) \geq 0\},
$$
is a non-empty subset of $\Omega$ defined by a single linear inequality.
Thus, each non-empty fiber is an intersection of a convex cone with
$\Sphere^s$.
Thus, all such fibres can be retracted to their 
respective centers of mass continuously,
proving the first half of the proposition.
\end{proof}

For any  quadratic form $Q$, we will denote by ${\rm index}(Q)$, the number of
negative eigenvalues of the symmetric matrix of the corresponding bilinear
form, that is of the matrix $M$ such that,
$Q(x) = \langle Mx, x \rangle$ for all $x \in \R^{k+1}$ 
(here $\langle\cdot,\cdot\rangle$ denotes the usual inner product). 
We will also
denote by $\lambda_i(Q), 0 \leq i \leq k$, the eigenvalues of $Q$, in non-decreasing order, i.e.
\[ \lambda_0(Q) \leq \lambda_1(Q) \leq \cdots \leq \lambda_k(Q).
\]

Given a 
quadratic map $P = (P_1,\ldots,P_s): \R^{k+1} \rightarrow \R^s,$
and $0 \leq j \leq k$, 
we denote by 
\[
\Omega_j = \{\omega \in \Omega \;  \mid \;  \lambda_j({\omega}P) \geq 0 \}.
\]

It is clear that the $\Omega_j$'s induce a filtration of the space
$\Omega$, i.e.,
$
\Omega_0 \subset \Omega_1 \subset \cdots \subset \Omega_{k}.
$

Agrachev \cite{Agrachev} showed that the Leray spectral sequence of the map
$\phi_1$ (converging to the cohomology $H^*(B) \cong H^*(T)$),
has as its $E_2$ terms,
\begin{equation}
\label{eqn:agrachev}
E_2^{pq} = H^p(\Omega_{k-q},\Omega_{k-q-1}).
\end{equation}
(We refer  the reader to \cite{Bredon} for the precise
definition of the Leray spectral sequence of a map). 

\label{obs:sphere}
Equation \ref{eqn:agrachev} follows directly from the following lemma.
\begin{lemma}
\label{lem:sphere}
The fibre of the map $\phi_1$ over a point 
$\omega \in \Omega_{j}\setminus \Omega_{j-1}$ has the homotopy type
of a sphere of dimension $k-j$. 
\end{lemma}

\begin{proof}
First notice that for
$\omega \in  \Omega_{j}\setminus \Omega_{j-1}$,
$\lambda_0({\omega}P),\ldots, \lambda_{j-1}({\omega}P) < 0.$ Moreover,
letting $Y_0({\omega}P),\ldots,Y_k({\omega}P)$ be 
the co-ordinates with respect to 
an orthonormal basis consisting of the
eigenvectors of ${\omega}P$, we have that 
$\phi_1^{-1}(\omega)$ is the subset of $\Sphere^k$ defined by,
$$
\displaylines{
\sum_{i=0}^{k} \lambda_i({\omega}P)Y_i({\omega}P)^2 \geq  0, \cr
\sum_{i=0}^{k} Y_i({\omega}P)^2 = 1.
}
$$

Since, $\lambda_i({\omega}P) < 0, 0 \leq i < j,$ it follows that
for $\omega \in \Omega_{j}\setminus \Omega_{j-1}$,
the fiber $\phi_1^{-1}(\omega)$ is homotopy equivalent to the
$(k-j)$-dimensional sphere defined by setting
$Y_0({\omega}P) = \cdots = Y_{j-1}({\omega}P) = 0$ on the sphere defined by
$\sum_{i=0}^{k}Y_i({\omega}P)^2 = 1.$
\end{proof}


\section{Computing the cohomology groups of a basic semi-algebraic
set defined by homogeneous quadratic inequalities
}
\label{sec:comp}
In this section, 
we will show how to effectively compute the spectral
sequence described in the previous section. 

Let ${\mathcal P}= (P_1,\ldots,P_s) \subset \R[X_0,\ldots,X_k]$ be a 
$s$-tuple of quadratic forms.
For any subset ${\mathcal Q} \subset {\mathcal P}$, we denote by
$T_{\mathcal Q} \subset \Sphere^k$, the semi-algebraic set,
$$
\displaylines{
T_{\mathcal Q} = \bigcup_{P \in {\mathcal Q}}
               \{x \in \Sphere^k\; \mid \; P(x) \leq 0 \},
}
$$
and let
$$
\displaylines{
S = \bigcap_{P \in {\mathcal P}}
               \{x \in \Sphere^k\; \mid \; P(x) \leq 0 \}.
}
$$
\noindent
We denote by $C^\bullet({\mathcal H}(T_{\mathcal Q}))$
the co-chain complex of a cellular subdivision, 
${\mathcal H}(T_{\mathcal Q})$ of $T_{\mathcal Q},$
which is to be chosen sufficiently fine  (to be specified later).

We first describe for each subset ${\mathcal Q} \subset {\mathcal P}$
with $\#{\mathcal Q} = \ell < k$,
a complex, ${\mathcal M}^{\bullet}_{\mathcal Q}$, 
and natural homomorphisms,
$$
\psi_{\mathcal Q}: C^\bullet({\mathcal H}(T_{\mathcal Q})) \rightarrow
{\mathcal M}^\bullet_{\mathcal Q},
$$
which induce isomorphisms,
$$
\psi_{\mathcal Q}^*:
H^{*}(C^{\bullet}({\mathcal H}(T_{\mathcal Q}))) \rightarrow
H^*({\mathcal M}_{\mathcal Q}^\bullet).
$$ 

Moreover, for ${\mathcal B} \subset {\mathcal A} \subset {\mathcal P}$ with
$\#{\mathcal A} = \#{\mathcal B} + 1  < k$, we construct a homomorphism
of complexes,
$$
\phi_{{\mathcal A},{\mathcal B}}: {\mathcal M}^{\bullet}_{\mathcal A} \rightarrow 
{\mathcal M}^{\bullet}_{\mathcal B},
$$
such that the following diagram commutes,

\begin{equation}
\label{eqn:commutative}
\begin{diagram}
H^*({\mathcal M}^{\bullet}_{\mathcal A})& \rTo^{\phi_{{\mathcal A},{\mathcal B}}^*}&
H^*({\mathcal M}^{\bullet}_{\mathcal B}) \\
\uTo^{\psi_{\mathcal A}^*}&&\uTo^{\psi_{\mathcal B}^*}\\
H^*(C^\bullet({\mathcal H}(T_{\mathcal A}))) &\rTo^{r^*} &H^*(C^\bullet({\mathcal H}(T_{\mathcal B})))
\end{diagram}
\end{equation}
\noindent
where $\phi_{{\mathcal A},{\mathcal B}}^*$ and $r^*$ are the induced homomorphisms
of $\phi_{{\mathcal A},{\mathcal B}}$ and the restriction homomorphism $r$ 
respectively.

Now, consider a fixed subset ${\mathcal Q} \subset {\mathcal P}$, which 
without loss of generality we take to be $\{P_1,\ldots,P_\ell\}$. Let
$$
P_{\mathcal Q} =(P_1,\ldots,P_\ell): \R^{k+1} \rightarrow \R^\ell
$$ 
denote the corresponding quadratic map.

As in the previous section, let 
$\R^{\mathcal Q} = \R^{\ell},$ and
$$
\Omega_{\mathcal Q} = \{\omega \in \R^{\ell} \mid  |\omega| = 1, 
\omega_i \leq 0, 1 \leq i \leq \ell\}.
$$

Let 
$B_{\mathcal Q} \subset \Omega_{\mathcal Q} \times \Sphere^k$ be the set defined by,
\[
B_{\mathcal Q} = \{ (\omega,x) \mid \omega \in \Omega_{\mathcal Q}, x \in \Sphere^k \;\mbox{and} \; 
{\omega}P_{\mathcal Q}(x) \geq  0\},
\]
and we denote by $\phi_{1,{\mathcal Q}}: {B}_{\mathcal Q} \rightarrow \Omega_{\mathcal Q}$ 
and 
$\phi_{2,{\mathcal Q}}: {B}_{\mathcal Q} \rightarrow \Sphere^k$ the two projection maps.

For each subset ${\mathcal Q}' \subset {\mathcal Q}$ we have a natural inclusion
$\Omega_{{\mathcal Q}'} \hookrightarrow \Omega_{\mathcal Q}$.

\subsection{Index Invariant Triangulations}
We now define a certain special kind of 
semi-algebraic triangulation of $\Omega_{\mathcal Q}$
that will play an important role later.
\begin{definition}(Index Invariant Triangulation)
\label{def:iit}
An {\em index invariant triangulation} of $\Omega_{\mathcal Q}$ consists of:
\begin{enumerate}
\item
A semi-algebraic triangulation,
$$
h_{\mathcal Q}: \Delta_{\mathcal Q} \rightarrow \Omega_{\mathcal Q}
$$ 
of  $\Omega_{\mathcal Q}$, 
which is compatible with the subsets $\Omega_{{\mathcal Q}'}$ for every
${\mathcal Q}' \subset {\mathcal Q}$, and
such that for any simplex $\sigma$ of $\Delta_{\mathcal Q}$, 
${\rm index}(\omega P_{\mathcal Q})$, 
as well as the multiplicities of the eigenvalues of $\omega P_{\mathcal Q}$,
stay invariant as $\omega$ varies over $h_{\mathcal Q}(\sigma)$, and
\item
for every simplex $\sigma$ of $\Delta_{\mathcal Q}$,
with ${\rm index}(\omega P_{\mathcal Q}) =  i$, 
a continuous choice of an orthonormal basis of $\R^{k+1}$,
$\{e_0(\sigma,\omega),\ldots,e_k(\sigma,\omega)\}$,
for $\omega \in h_{\mathcal Q}(\sigma)$,
such that $\{e_i(\sigma,\omega),\ldots,e_k(\sigma,\omega)\}$ span 
the linear subspace of $\R^{k+1}$ on which the quadratic
form $\omega P_{\mathcal Q}$ is positive semi-definite.
\end{enumerate}
\end{definition}

We describe later how 
to compute index invariant triangulations
(see Algorithm \ref{alg:triangulation} below). 
It will follow from the complexity analysis of 
Algorithm \ref{alg:triangulation} that
the size of the complex $\Delta_{\mathcal Q}$ as well as the
degrees of the polynomials occuring in the parametrized univariate representations
defining the various bases,
$\{e_0(\sigma,\omega),\ldots,e_k(\sigma,\omega)\}$,
are all bounded by $k^{2^{O(\ell)}}$.

For the rest of this section we fix an index invariant triangulation
$
h_{\mathcal Q}: \Delta_{\mathcal Q} \rightarrow \Omega_{\mathcal Q},
$ 
satisfying the complexity estimates stated above.

The following proposition relates,
for any simplex $\sigma \in \Delta_{\mathcal Q}$,
the homotopy type of 
$\phi_{1,{\mathcal Q}}^{-1}(h_{\mathcal Q}(\sigma))$ to that of a single fiber.
 \begin{proposition}
\label{prop:homotopy1}
For any simplex $\sigma \in \Delta_{\mathcal Q}$ and 
$\omega \in h_{\mathcal Q}(\sigma),$
$\phi_{1,{\mathcal Q}}^{-1}(h_{\mathcal Q}(\sigma))$ is homotopy equivalent to
$\phi_{1,{\mathcal Q}}^{-1}(\omega),$ 
and both these spaces have the homotopy type of the 
sphere $\Sphere^{k - {\rm index}(\omega P)}.$ 
\end{proposition}

\begin{proof}
Let $i = {\rm index}(\omega  P).$
Since ${\rm index}(\omega P)$ is invariant as $\omega$ varies over
$h_{\mathcal Q}(\sigma)$, the quadratic forms $\omega P$ has exactly $i$
negative eigen-values for each $\omega \in h_{\mathcal Q}(\sigma).$
Let $M(\sigma,\omega) \subset \R^{k+1}$ be 
the 
orthogonal complement to the 
linear span of the corresponding eigen-vectors,
and let $B(\sigma,\omega) = M(\sigma,\omega) \cap \Sphere^k.$
Clearly, $M(\sigma,\omega)$ and $B(\sigma,\omega)$ 
vary continuously with $\omega$,
and $\phi_{1,{\mathcal Q}}^{-1}(\omega)$ can be retracted to the set 
$\{\omega\} \times B(\sigma,\omega).$
Finally, since $h_{\mathcal Q}(\sigma)$ is contractible to $\omega$,
its clear that $\phi_{1,{\mathcal Q}}^{-1}(h_{\mathcal Q}(\sigma))$ retracts to
$\{\omega\} \times  B(\sigma,\omega)$ and the latter has the
homotopy type of $\Sphere^{k - {\rm index}(\omega  P)}$ by 
Lemma \ref{lem:sphere}. 
\end{proof}

\subsection{Definition of the cell complex ${\mathcal K}(B_{\mathcal Q})$}
\label{subsec:cellcomplex}
Our next goal is to construct a cell complex homotopy equivalent to
$B_{\mathcal Q}$ obtained by glueing together certain regular 
cell complexes,  ${\mathcal K}(\sigma)$,
where $\sigma \in \Delta_{\mathcal Q}$.

       \begin{figure}[hbt]
         \centerline{
           \scalebox{0.5}{
             \input{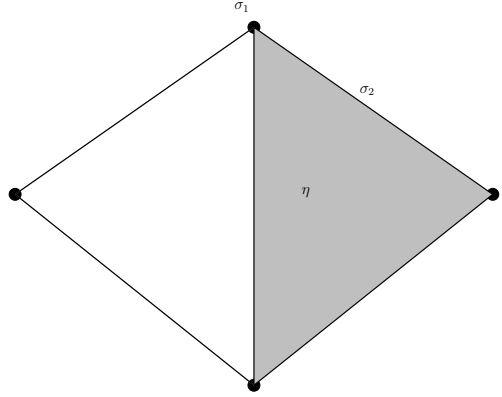}
             }
           }
         \caption{The complex $\Delta_{\mathcal Q}$.}
         \label{fig-eg1}
       \end{figure}

       \begin{figure}[hbt]
         \centerline{
           \scalebox{0.5}{
             \input{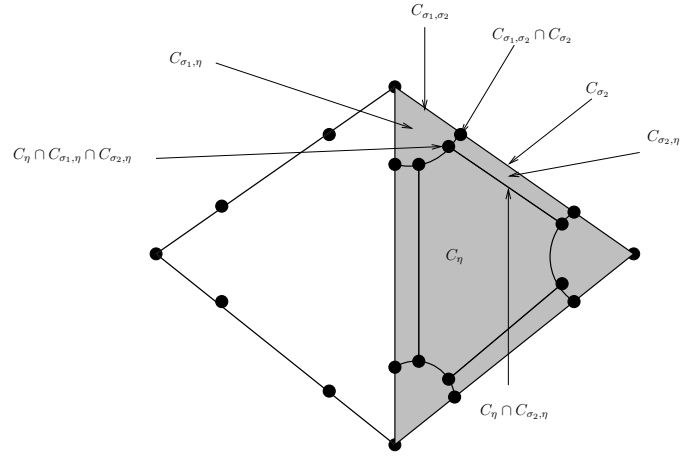}
             }
           }
         \caption{The corresponding complex 
${\mathcal C}(\Delta_{\mathcal Q})$.}
         \label{fig-eg2}
       \end{figure}

Let $1 \geq  \eps_0 \gg \eps_1 \gg \cdots \gg \eps_s \gg 0$ be 
infinitesimals. 
For $\eta \in \Delta_{\mathcal Q}$, 
we denote by $C_{\eta}$ the subset of $\bar\eta$ defined by,
\[
C_{\eta} = \{x \in \bar\eta \;\mid\; 
\mbox{ and }\dist(x,\theta) \geq \eps_{\dim(\theta)} \mbox{ for all }
\theta \prec \sigma \}.
\]

Now, let $\sigma \prec \eta$ be two simplices of $\Delta_{\mathcal Q}$.
We denote by 
$C_{\sigma,\eta}$ the subset of $\bar\eta$ defined by,
\[
C_{\sigma,\eta} = \{x \in \bar\eta \;\mid\; \dist(x,\sigma) \leq 
\eps_{\dim(\sigma)},
\mbox{ and }\dist(x,\theta) \geq \eps_{\dim(\theta)} \mbox{ for all }
\theta \prec \sigma \}.
\]

Note that,
$$
\displaylines{
|\Delta_{\mathcal Q}| = 
\bigcup_{\sigma \in \Delta_{\mathcal Q}} C_{\sigma} \cup
\bigcup_{\sigma,\eta \in \Delta_{\mathcal Q},\sigma \prec \eta} C_{\sigma,\eta}.
}
$$

Also, observe that the various $C_\eta$'s and $C_{\sigma,\eta}$'s 
are all homeomorphic to balls,
and moreover all non-empty intersections between them also have the same property.
Thus, the union of the $C_{\eta}$'s and $C_{\sigma,\eta}$'s together with the
non-empty intersections between them form  a regular cell complex,
${\mathcal C}(\Delta_{\mathcal Q})$, whose underlying
topological space is $|\Delta_{\mathcal Q}|$ 
(see Figures \ref{fig-eg1} and \ref{fig-eg2}).

We now associate to each 
$C_{\sigma}$  (respectively,  $C_{\sigma,\eta}$)
a regular cell complex, ${\mathcal K}(\sigma)$, (respectively,
${\mathcal K}(\sigma,\eta)$)
homotopy equivalent to 
$\phi_{1,{\mathcal Q}}^{-1}(h_{\mathcal Q}(C_\sigma))$
(respectively, \\
$
\displaystyle{
\phi_{1,{\mathcal Q}}^{-1}(h_{\mathcal Q}(C_{\sigma,\eta})).
}
$

For each $\sigma \in \Delta_{\mathcal Q}$, and 
$\omega \in h_{\mathcal Q}(\sigma)$, let 
$$
\displaylines{
\lambda_0^{\sigma}(\omega) = \cdots = \lambda_{i_0}^{\sigma}(\omega)
< \lambda_{i_0+1}^{\sigma}(\omega) = \cdots = \lambda_{i_1}^\sigma(\omega)
< \cdots < \lambda_{i_{p-1}+1}^\sigma(\omega) = \cdots = \cr
\lambda_{i_p}^{\sigma}(\omega) < 0 = \cdots = 
\lambda_{i_{p+1}}^{\sigma}(\omega) < \cdots = \cdots = 
\lambda_k^{\sigma}(\omega)
}
$$
denote the eigenvalues of $\omega P$. Here, ${\rm index}(\omega P) = i_p +  1$.
Also, since the multiplicities of the eigenvalues do not change as
$\omega$ varies over $h_{\mathcal Q}(\sigma)$, the block structure,
$[0,\ldots,i_0], [i_0+1,\ldots,i_1],\ldots,[\cdot,\ldots,k]$ also does
not change as $\omega$ varies over $h_{\mathcal Q}(\sigma)$.
For $0 \leq j \leq p$, let $M^j(\sigma,\omega)$ denote the subspace of
$\R^{k+1}$ orthogonal to the subspace spanned by the eigenvectors
corresponding to the eigenvalues 
$\lambda_0^{\sigma}(\omega) = \cdots = \lambda_{i_0}^{\sigma}(\omega)
< \lambda_{i_0+1}^{\sigma}(\omega) = \cdots = \lambda_{i_1}^\sigma(\omega)
< \cdots < \lambda_{i_{p-1}+1}^\sigma(\omega) = \cdots = 
\lambda_{i_j}^{\sigma}(\omega)$,
and let  $M(\sigma,\omega) = M^{p}(\sigma,\omega)$. 
Since the eigenvalues vary continuously
and their multiplicities do not change as $\omega$ varies over $h_{\mathcal Q}(\sigma)$,
the flag of subspaces $M^0(\sigma,\omega)  \supset \cdots \supset 
M^p(\sigma,\omega)$ also
varies continuously over $h_{\mathcal Q}(\sigma)$. 

For each $\sigma \in \Delta_{\mathcal Q}$, 
and $\omega \in h_{\mathcal Q}(\sigma)$, let
$\{e_0(\sigma,\omega),\ldots,e_{k}(\sigma,\omega)\},$ be 
the continuously varying  orthonormal basis of $\R^{k+1}$
computed previously.

We extend the the orthonormal basis,
$\{e_0(\sigma,\omega),\ldots,e_{k}(\sigma,\omega)\},$ 
continuously to
each $C_{\sigma,\eta}$ for $\eta$ with $\sigma \prec \eta$,
satisfying the condition that 
\[
M(\sigma,\omega)  \subset 
\spanof(e_i(\sigma,\omega),\ldots,e_{k}(\sigma,\omega)).
\] 
This extension can be done in a consistent manner because
the first $i$ eigenvalues, $\lambda_0(\omega P),
\ldots, \lambda_{i-1}(\omega P)$  of $\omega P$ stay negative,
and $\lambda_{i-1}(\omega P) < \lambda_{i}(\omega P)$
for $\omega$ in any infinitesimal neighborhood of 
$h_{\mathcal Q}(C_{\sigma})$. Thus, the linear subspace of $\R^k$
orthogonal to the eigenspaces corresponding to the eigenvalues, 
$\lambda_0(\omega P), \ldots, \lambda_{i-1}(\omega P)$ is
well defined, 
and varies continuously with  $\omega$ in any infinitesimal neighborhood of 
$h_{\mathcal Q}(C_{\sigma})$.

More precisely, for any point $z' = h_{\mathcal Q}^{-1}(\omega')
\in C_{\sigma,\eta}$, let 
$z = h_{\mathcal Q}^{-1}(\omega)$, 
be the unique point in $C_{\sigma}$ closest to $z'$.
For $j = i,\ldots,k$, let $e_j'(\sigma,\omega')$ be the orthogonal
projection of $e_j(\sigma,\omega)$ onto the subspace 
$M(\sigma,\omega')$, and let $e_i(\sigma,\omega'),\ldots,e_k(\sigma,\omega')$
be obtained from $e_i'(\sigma,\omega'),\ldots,e_k'(\sigma,\omega')$ by
Gram-Schmidt orthogonalization. 
Note that each $e_j(\sigma,\omega')$ obtained this way is infinitesimally
close to $e_j(\sigma,\omega)$.
Using the same procedure starting with the
vectors $e_0(\sigma,\omega),\ldots, e_{i-1}(\sigma,\omega)$ and the
linear subspace $M(\sigma,\omega')^{\perp}$, we compute
an orthonormal set of vectors  
$e_0(\sigma,\omega'),\ldots,e_{i-1}(\sigma,\omega')$ spanning
$M(\sigma,\omega')^{\perp}$. It is clear from construction, that
the vectors $e_0(\sigma,\omega'),\ldots,e_k(\sigma,\omega')$ form
an orthonormal basis for $\R^{k+1}$. Moreover, they
vary continuously with $\omega'$, extending the basis
$e_0(\sigma,\omega),\ldots,e_k(\sigma,\omega)$, and finally
$e_i(\sigma,\omega'),\ldots,e_k(\sigma,\omega')$ span
$M(\sigma,\omega')$.

The orthonormal basis 
\[
\{e_0(\sigma,\omega),\ldots,e_{k}(\sigma,\omega)\},
\] 
determines a complete flag of subspaces, 
${\mathcal F}(\sigma,\omega)$, consisting of
$$
\displaylines{
L^0(\sigma,\omega) = 0, \cr
L^1(\sigma,\omega) = \spanof(e_k(\sigma,\omega)),\cr
 L^2(\sigma,\omega) = 
\spanof(e_k(\sigma,\omega),e_{k-1}(\sigma,\omega)), \cr
\vdots \cr
L^{k+1}(\sigma,\omega) = \R^{k+1}.
}
$$

For $0 \leq j \leq k$, let $c_{j}^+(\sigma,\omega)$ 
(respectively, $c_{j}^-(\sigma,\omega)$)
denote the $(k-j)$-dimensional cell consisting of the intersection of the
$L^{k-j+1}(\sigma,\omega)$
with the unit hemisphere in $\R^{k+1}$ 
defined by $\{x \in \Sphere^k\mid \langle x,e_j(\sigma,\omega)\rangle \geq 0\}$
(respectively, $\{x \in \Sphere^k\mid \langle x,e_j(\sigma,\omega)
\rangle \leq 0\}$).

The regular cell complex ${\mathcal K}(\sigma)$ 
(as well as ${\mathcal K}(\sigma,\eta)$)
is defined as follows.

The cells of ${\mathcal K}(\sigma)$ are
$\{(x,\omega) \mid x \in c_j^{\pm}(\sigma,\omega), \omega \in h_{\mathcal Q}(c)\}$,
where ${\rm index}(\omega P) \leq j \leq k$,
and 
$c \in {\mathcal C}(\Delta_{\mathcal Q})$
is either $C_\sigma$ itself, or a cell
contained in the boundary  of $C_\sigma$.

Similarly, the cells of ${\mathcal K}(\sigma,\eta)$ are
$\{(x,\omega) \mid x \in c_j^{\pm}(\sigma,\omega), \omega \in h_{\mathcal Q}(c)\}$,
where 
${\rm index}(\omega P) \leq j \leq k$,
$c \in {\mathcal C}(\Delta_{\mathcal Q})$
is either $C_{\sigma,\eta}$ itself,  or a cell
contained in the boundary  of $C_{\sigma,\eta}$.

Our next step is to obtain cellular subdivisions
of each non-empty intersection amongst the
spaces associated to the complexes constructed above, and thus obtain
a regular cell complex,
${\mathcal K}(B_{\mathcal Q})$, whose associated space,
$|{\mathcal K}(B_{\mathcal Q})|$, will be shown to be 
homotopy equivalent to $B_{\mathcal Q}$ (Proposition \ref{prop:iso2} below).

First notice that $|{\mathcal K}(\sigma',\eta')|$ (respectively, 
$|{\mathcal K}(\sigma)|$) has a non-empty intersection with 
$|{\mathcal K}(\sigma,\eta)|$ only if $C_{\sigma',\eta'}$ (respectively,
$C_{\sigma'}$) intersects $C_{\sigma,\eta}$. 

Let $C$ be some non-empty intersection amongst the 
$C_{\sigma}$'s and $C_{\sigma,\eta}$'s,
that is $C$ is a cell of ${\mathcal C}(\Delta_{\mathcal Q})$.
Then, $C \subset \eta$ for a unique simplex $\eta \in \Delta_{\mathcal Q}$, 
and 
$$
\displaylines{
C =  C_{\sigma_1,\eta} \cap \cdots \cap C_{\sigma_p,\eta},
}
$$
with $\sigma_1 \prec \sigma_2 \prec \cdots \prec \sigma_p \prec \eta$
and $p \leq \#{\mathcal Q}+1$.

Consider $\omega \in h_{\mathcal Q}(C)$. 
We have $p$ different flags, 
\[
{\mathcal F}(\sigma_1,\omega),
\ldots, {\mathcal F}(\sigma_p,\omega),
\]
giving rise to  $p$ independent regular cell decompositions of 
$
\displaystyle{
B(\omega,\eta) = M(\eta,\omega) \cap \Sphere^{k}.
}
$

       \begin{figure}[h]
         \centerline{
           \scalebox{0.5}{
             \input{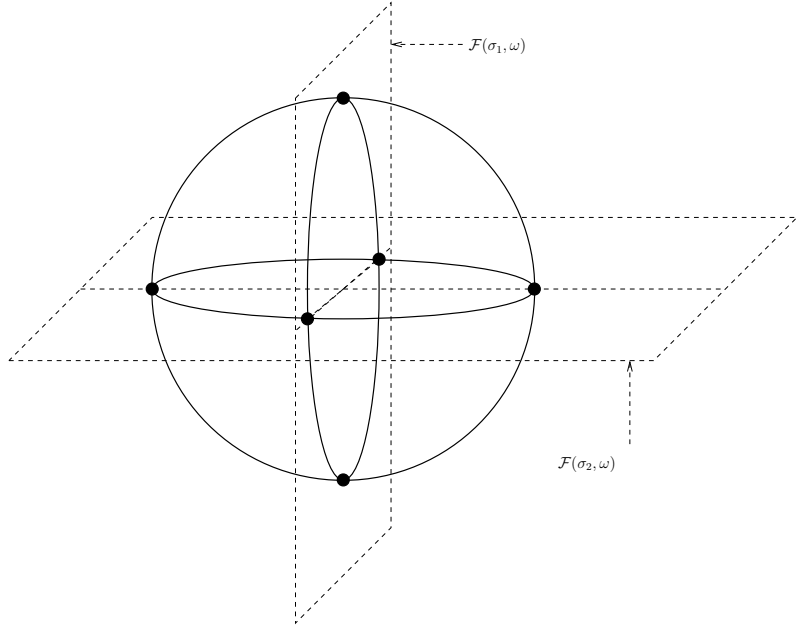}
             }
           }
         \caption{The cell complex ${\mathcal K}'(C,\omega)$.}  
         \label{fig-eg5}
       \end{figure}

There is a unique smallest regular cell complex, ${\mathcal K}'(C,\omega)$,  
that refines all these cell decompositions.
The cells of this cell decomposition consists of the following.
Let $L \subset M(\eta,\omega)$ be  any linear subspace of dimension 
$m, 0 \leq m \leq k+1$, which is an intersection,
of linear subspaces $L_1,\ldots,L_p$, where 
$L_i \in {\mathcal F}(\sigma_i,\omega), 1 \leq i \leq p$. 
The elements of the flags,
${\mathcal F}(\sigma_1,\omega),
\ldots, {\mathcal F}(\sigma_p,\omega)$ of dimensions $m+1$,
partition $L$ into polyhedral cones of various dimensions. The
union of the sets of intersections of these cones with $\Sphere^{k}$, 
over all such subspaces $L \subset M(\eta,\omega)$, are the cells of 
${\mathcal K}'(C,\omega)$.
Figure \ref{fig-eg5} illustrates the refinement described above in case
of two flags in $\R^3$.

We now triangulate $h_{\mathcal Q}(C)$, using 
the algorithm implicit in Theorem  \ref{the:triangulation} (Triangulation),
such that the combinatorial type of the arrangement of flags, 
\[
{\mathcal F}(\sigma_1,\omega),
\ldots, {\mathcal F}(\sigma_p,\omega)
\] 
and hence the cell decomposition ${\mathcal K}'(C,\omega)$,  
stays invariant over the image,
$h_C(\theta)$,  of  each simplex, $\theta$, of this triangulation. 
More precisely, we first compute a family of polynomials,
${\mathcal A}_{C}\subset \R[Z_1,\ldots,Z_\ell]$ 
whose signs at $\omega$ determine the combinatorial type of the
corresponding arrangement of flags. It is easy to verify
(see Section \ref{subsec:arrangement}),
given the complexity bounds on the parametrized univariate 
representations defining the orthonormal bases,
$\{e_0(\sigma,\omega),\ldots,e_k(\sigma,\omega)\}$,
$\omega \in h_{\mathcal Q}(\sigma)$, stated above,
that the
number and degrees of the polynomials in the family ${\mathcal A}_C$
is bounded by $k^{2^{O(\ell)}}$.
We then
use the algorithm implicit in Theorem  \ref{the:triangulation} (Triangulation),
with ${\mathcal A}_{C}$ as input, to obtain the required triangulation.

The closures of the sets
\[
\{(\omega,x) \;\mid\; x \in c \in  {\mathcal K}'(C,\omega),
\omega \in h_{\mathcal Q}(h_C(\theta))\}
\] 
constitute a regular cell complex, 
${\mathcal K}(C)$, which is compatible with the regular cell complexes
${\mathcal K}(\sigma_1),\ldots, {\mathcal K}(\sigma_p)$.

The following proposition gives an upper bound on the size of the
complex ${\mathcal K}(C)$. We use the notation introduced in the previous
paragraph.
\begin{proposition}
\label{prop:complexity}
For each $\omega \in h_{\mathcal Q}(C)$, the number of cells in ${\mathcal K}'(C,\omega)$
is bounded by $k^{O(\ell)}$. Moreover, the number of cells in the complex
${\mathcal K}(C)$ is bounded by $k^{2^{O(\ell)}}$.
\end{proposition}

\begin{proof}
The first part of the proposition follows from the fact that there are at
most $k^{\#{\mathcal Q}+1} = k^{\ell+1}$ 
choices for the linear space $L$ and the number of $(m-1)$
dimensional cells contained in $L$ is bounded by $2^{\ell}$ (which is an upper
bound on the number of full dimensional cells in an arrangement of at most
$\ell$ hyperplanes).
The second part is a consequence of the complexity estimate
in Theorem \ref{the:triangulation} (Triangulation) and the bounds on
number and degrees of polynomials in the family ${\mathcal A}_C$
stated above.
\end{proof}

We denote by ${\mathcal K}(B_{\mathcal Q})$, the union of all the complexes
${\mathcal K}(C)$ constructed above, noting that by construction,
${\mathcal K}(B_{\mathcal Q})$ is a regular cell complex.

\begin{proposition}
\label{prop:iso2}
$|{\mathcal K}(B_{\mathcal Q})|$ is homotopy equivalent to $B_{\mathcal Q}$.
\end{proposition}
\begin{proof}
We first show that
$B_{\mathcal Q}$ is homotopy equivalent to a subset 
$B_{\mathcal Q}' \subset B_{\mathcal Q}$ as follows.
For each simplex $\sigma \in \Delta_{\mathcal Q}$ of the largest dimension
$\ell$, we use the retraction used in the proof of 
Proposition \ref{prop:homotopy1}, to retract 
$\phi_1^{-1}(h_{\mathcal Q}({\rm rel int}(C_{\sigma})))$ to the set 
$\{(\omega,x) \;\mid\; \omega \in {\rm rel int}(C_{\sigma}), 
x \in B(\omega,\sigma)\}$.
In this way we obtain a semi-algebraic set, $X_\ell'$, 
which is a deformation
retract of $\E(B_{\mathcal Q},\R\la\eps_0,\ldots,
\eps_{\ell-1}\ra)$. 
Let  $X_\ell = \lim_{\eps_{\ell-1}} X_\ell'$. 
Notice that in the definition of $X_\ell'$, if we
replace $\eps_{\ell-1}$ by a variable $t$ and denote the corresponding set by
$X_{\ell,t}'$, then for all $0 < t < t'$, $X_{\ell,t}' \subset X_{\ell,t'}'$
and each $X_{\ell,t}$ is closed and bounded.
It then follows (see Lemma 16.17 in \cite{BPR03}) that
$\E(X_\ell,\R\la\eps_0,\ldots,
\eps_{\ell-1}\ra)$ 
has the same homotopy type as $X_{\ell}'$, and hence 
$X_\ell$ has the same homotopy type as $\E(B_{\mathcal Q},\R\la\eps_0,\ldots,
\eps_{\ell-2}\ra)$. 

Now repeat the process using the $({\ell}-1)$-dimensional 
simplices and so on, to finally obtain
$X_0 = B_{\mathcal Q}'$, which by construction has the same homotopy type
as $B_{\mathcal Q}$.
Finally, (again using Lemma 16.17 in \cite{BPR03}) we also have that
$X_0 = \lim_{\eps_0} |{\mathcal K}(B_{\mathcal Q})|$ and 
$\E(X_0, \R\la\eps_0,\ldots,\eps_{\ell-1}\ra)$ 
has the same homotopy type as $|{\mathcal K}(B_{\mathcal Q})|$.
\end{proof}

We also have,

\begin{proposition}
\label{prop:complexity2}
The number of cells in the
cell complex ${\mathcal K}(B_{\mathcal Q})$ is 
bounded by $k^{2^{O(\ell)}}$.
\end{proposition}
\begin{proof}
The proposition is a consequence of Proposition \ref{prop:complexity} and
the fact that the number of cells in the complex 
${\mathcal C}(\Delta_{\mathcal Q})$ is bounded by $k^{2^{O(\ell)}}.$
\end{proof}

We now define,
\[
{\mathcal M}_{\mathcal Q}^{\bullet}
= C^{\bullet}({\mathcal K}(B_{\mathcal Q}),
\]
where 
$C^{\bullet}({\mathcal K}(B_{\mathcal Q})$ is the cellular co-chain
complex of the regular cell complex ${\mathcal K}(B_{\mathcal Q})$.

Let ${\mathcal H}(T_{\mathcal Q})$ 
(resp. ${\mathcal H}(B_{\mathcal Q})$)
be a suitably fine cellular subdivision
of $T_{\mathcal Q}$ 
(resp. $B_{\mathcal Q}$)
and let 
\[
\phi_{2,{\mathcal Q}}': 
C_\bullet({\mathcal H}({B}_{\mathcal Q}))
\rightarrow 
C_\bullet({\mathcal H}(T_{\mathcal Q})),
\]
be the homomorphism induced by a cellular map, 
which is a cellular approximation of $\phi_{2,{\mathcal Q}}$.

Let
$\phi_{\mathcal Q}:  |{\mathcal K}(B_{\mathcal Q})| \rightarrow 
B_{\mathcal Q}$ denote the  homotopy equivalence 
shown to exist by Proposition \ref{prop:iso2} above and
let
\[
\phi_{\mathcal Q}': 
C_\bullet({\mathcal K}'(B_{\mathcal Q})) 
\rightarrow 
C_\bullet({\mathcal H}(B_{\mathcal Q})),
\]
be the homomorphism induced by a cellular approximation to $\phi_{\mathcal Q}$,
where ${\mathcal K}'(B_{\mathcal Q})$ is a cellular subdivision of the
complex ${\mathcal K}(B_{\mathcal Q})$.

Since, each cell of ${\mathcal K}(B_{\mathcal Q})$ is a union of cells
of ${\mathcal K}'(B_{\mathcal Q})$, there is a natural homomorphism
\[
\theta_{\mathcal Q}: C_\bullet({\mathcal K}(B_{\mathcal Q})) 
\rightarrow C_\bullet({\mathcal K}'(B_{\mathcal Q})) 
\]
obtained by sending each $p$-dimensional cell of 
${\mathcal K}(B_{\mathcal Q})$ to the sum of $p$-dimensional cells
of  ${\mathcal K}'(B_{\mathcal Q})$ contained in it, for every $p \geq 0$.
It is a standard fact that $\theta_{\mathcal Q}$ and its dual,
$\check{\theta}_{\mathcal Q}$, are quasi-isomorphisms.

Let 
\[
\psi_{\mathcal Q} = \check{\theta}_{\mathcal Q}\circ  \check{\phi}_{\mathcal Q}'\circ \check{\phi}_{2,{\mathcal Q}}': 
C^\bullet({\mathcal H}(T_{\mathcal Q})) \rightarrow 
C^\bullet({\mathcal K}(B_{\mathcal Q})),
\]
where $\check{\phi}_{\mathcal Q}'$ (resp. $ \check{\phi}_{2,{\mathcal Q}}'$)
is the dual homomorphism  of 
${\phi}_{\mathcal Q}'$ (resp. $ {\phi}_{2,{\mathcal Q}}'$).

\begin{proposition}
\label{prop:iso}
For $0 \leq i \leq k-1,$ the induced homomorphisms,
$$
\psi_{\mathcal Q}^*: H^i(C^{\bullet}({\mathcal H}(T_{\mathcal Q}))) 
\rightarrow H^i({\mathcal M}^\bullet_{\mathcal Q})
$$
are isomorphisms.
\end{proposition}
\begin{proof}
The proof is clear since $\psi_{\mathcal Q}$ is a composition of
quasi-isomorphisms.
\end{proof}

Now let,
${\mathcal B} \subset {\mathcal A} \subset {\mathcal P}$ with
$\#{\mathcal A} = \#{\mathcal B} + 1  < k$. 

The simplicial complex $\Delta_{\mathcal B}$ is a subcomplex of 
$\Delta_{\mathcal A}$ and hence,
${\mathcal K}(B_{\mathcal B})$
is a subcomplex of 
${\mathcal K}(B_{\mathcal A})$ and thus there exists a natural
homomorphism (induced by restriction),
$$
\phi_{{\mathcal A},{\mathcal B}}: 
{\mathcal M}_{\mathcal A}^{\bullet}
\rightarrow
{\mathcal M}_{\mathcal B}^{\bullet}.
$$

The complexes ${\mathcal M}_{\mathcal A}^{\bullet},
{\mathcal M}_{\mathcal B}^{\bullet}$, and the homomorphisms, 
$
\phi_{{\mathcal A},{\mathcal B}},
\psi_{\mathcal A}, \psi_{\mathcal B}
$ 
satisfy
 
\begin{proposition}
\label{prop:commutative}
The diagram
\begin{diagram}
{\mathcal M}^{\bullet}_{\mathcal A}& \rTo^{\phi_{{\mathcal A},{\mathcal B}}}&
{\mathcal M}^{\bullet}_{\mathcal B} \\
\uTo^{\psi_{\mathcal A}}&&\uTo^{\psi_{\mathcal B}}\\
C^\bullet({\mathcal H}(T_{\mathcal A})) &\rTo^{r} & 
C^\bullet({\mathcal H}(T_{\mathcal B}))
\end{diagram}
is commutative,
where $r$ is the restriction homomorphism.
\end{proposition}

\begin{proof}
Clear from the construction.
\end{proof}

It follows from Proposition \ref{prop:commutative} that the diagram
(\ref{eqn:commutative}) is also commutative.

We denote by
$$ 
\check{\phi}_{{\mathcal B},{\mathcal A}}:
\check{{\mathcal M}}_{\mathcal B}^{\bullet}
\rightarrow
\check{{\mathcal M}}_{\mathcal A}^{\bullet}
$$
the homomorphism dual to $\phi_{{\mathcal A},{\mathcal B}}$.
We denote by ${\mathcal D}^{\bullet,\bullet}_{\mathcal P}$ 
the double complex defined by:

$$
\displaylines{
{\mathcal D}_{\mathcal P}^{p,q} = 
\bigoplus_{{\mathcal Q} \subset {\mathcal P},\#{\mathcal Q} = p+1}
{\check{\mathcal M}}^{q}_{\mathcal Q}.
}
$$

The vertical differentials,
$$
\displaylines{
d: {\mathcal D}_{\mathcal P}^{p,q} \rightarrow {\mathcal D}_{\mathcal P}^{p,q-1},
}
$$
are induced componentwise from the
differentials of the individual complexes 
${\check{\mathcal M}}^{\bullet}_{\mathcal Q}$.
The horizontal differentials,
$$
\displaylines{
\delta: {\mathcal D}_{\mathcal P}^{p,q} \rightarrow {\mathcal D}_{\mathcal P}^{p+1,q},
}
$$
are defined as follows:
for $a \in {\mathcal D}_{\mathcal P}^{p,q} = 
\oplus_{\#{\mathcal Q} = p+1}{\check{\mathcal M}}^q_{\mathcal Q},$
and for each subset 
\[
{\mathcal Q} = \{P_{i_0},\ldots,P_{i_{p+1}}\} \subset {\mathcal P}
\] 
with $i_0 < \cdots < i_{p+1},$ 
the ${\mathcal Q}$-th component of 
$\delta a \in {\mathcal D}_{\mathcal P}^{p+1,q}$ is given by,
$$
(\delta a)_{{\mathcal Q}} = \sum_{0 \leq j \leq p+1}
                         \check{\phi}_{{\mathcal Q}_j,{\mathcal Q}}(a_{{\mathcal Q}_j}),
$$

where ${\mathcal Q}_j = {\mathcal Q} \setminus \{P_{i_j}\}.$
{\small
$$
\begin{array}{ccccccccc}
& & \vdots  && \vdots  && \vdots  && \cr
& &
\Big\downarrow\vcenter{\rlap{$d$}} & &
\Big\downarrow\vcenter{\rlap{$d$}} & &
\Big\downarrow\vcenter{\rlap{$d$}} & & \cr
0 & \longrightarrow & \oplus_{\#{\mathcal Q} = 1}{\check{\mathcal M}}^{3}_{\mathcal Q} & 
\stackrel{\delta}{\longrightarrow}  & \oplus_{\#{\mathcal Q} = 2}{\check{\mathcal M}}^3_{\mathcal Q}&
\stackrel{\delta}{\longrightarrow} & \oplus_{\#{\mathcal Q} = 3}{\check{\mathcal M}}^3_{\mathcal Q}&
 \longrightarrow & \cdots
\cr
& &
\Big\downarrow\vcenter{\rlap{$d$}} & &
\Big\downarrow\vcenter{\rlap{$d$}} & &
\Big\downarrow\vcenter{\rlap{$d$}} & &\cr
0 & \longrightarrow & \oplus_{\#{\mathcal Q} = 1}{\check{\mathcal M}}^2_{\mathcal Q} &
\stackrel{\delta}{\longrightarrow} & \oplus_{\#{\mathcal Q} = 2}{\check{\mathcal M}}^2_{\mathcal Q} &
\stackrel{\delta}{\longrightarrow} & \oplus_{\#{\mathcal Q} = 3}{\check{\mathcal M}}^2_{\mathcal Q} &
 \longrightarrow &\cdots 
\cr
& &
\Big\downarrow\vcenter{\rlap{$d$}} & &
\Big\downarrow\vcenter{\rlap{$d$}} & &
\Big\downarrow\vcenter{\rlap{$d$}} & &\cr
0 & \longrightarrow & \oplus_{\#{\mathcal Q} = 1}{\check{\mathcal M}}^1_{\mathcal Q}&
\stackrel{\delta}{\longrightarrow} & \oplus_{\#{\mathcal Q} = 2}{\check{\mathcal M}}^1_{\mathcal Q}&
\stackrel{\delta}{\longrightarrow} & \oplus_{\#{\mathcal Q} = 3}{\check{\mathcal M}}^1_{\mathcal Q} &
 \longrightarrow & \cdots
\cr
& &
\Big\downarrow\vcenter{\rlap{$d$}} & &
\Big\downarrow\vcenter{\rlap{$d$}} & &
\Big\downarrow\vcenter{\rlap{$d$}} & & \cr
0 & \longrightarrow & \oplus_{\#{\mathcal Q} = 1}{\check{\mathcal M}}^0_{\mathcal Q} &
\stackrel{\delta}{\longrightarrow} & \oplus_{\#{\mathcal Q} = 2}{\check{\mathcal M}}^0_{\mathcal Q} &
\stackrel{\delta}{\longrightarrow} & \oplus_{\#{\mathcal Q} = 3}{\check{\mathcal M}}^0_{\mathcal Q} &
 \longrightarrow & \cdots
\cr
& &
\Big\downarrow\vcenter{\rlap{$d$}} & &
\Big\downarrow\vcenter{\rlap{$d$}} & &
\Big\downarrow\vcenter{\rlap{$d$}} & & \cr
& & 0 && 0 && 0 & &\cr 
\end{array}
$$
}

We have the following theorem.

\begin{theorem}
\label{the:main}
For $0 \leq i \leq k,$
$$
H^i(S) \cong H^i({\rm Tot}^{\bullet}({\mathcal D}^{\bullet,\bullet}_{\mathcal P})).
$$
\end{theorem}
\begin{proof}
For  each $i, 1 \leq i \leq s$, let $S_i \subset \Sphere^k$ denote,
the set defined on $\Sphere^k$ by $P_i \leq 0$. 
Then, $S = \cap_{i=1}^s S_i.$
Choosing a suitably fine triangulation of $\cup_{i=1}^s S_i$ consider the
Mayer-Vietoris double complex, ${\mathcal N}^{\bullet,\bullet}$, as described in
Section \ref{subsec:MV}.
The homomorphisms,
$$
\bigoplus_{{\mathcal Q} \subset {\mathcal P},
\#{\mathcal Q} = p+1} \check{\psi}_{\mathcal Q}: 
\bigoplus_{{\mathcal Q} \subset {\mathcal P},
\#{\mathcal Q} = p+1}
\check{\mathcal M}^q_{\mathcal Q}
\longrightarrow  \bigoplus_{{\mathcal Q} \subset {\mathcal P},
\#{\mathcal Q} = p+1} C_q({\mathcal H}(T_{\mathcal Q}))
$$
give a homomorphism of the double complexes,
\begin{equation}
\label{eqn:fundamental}
\psi: {\mathcal D}_{\mathcal P}^{\bullet,\bullet} 
\longrightarrow 
{\mathcal N}^{\bullet,\bullet}.
\end{equation}
By Proposition \ref{prop:iso}, $\psi$ induces isomporphisms
between the $E_1$ terms of the two spectral sequences, obtained
by taking homology with respect to the vertical differentials.
Theorem \ref{the:spectral} then implies that
$\psi$ induces isomorphisms 
in the associated spectral sequences. But, the second spectral
sequence converges to the homology of $S$ by (\ref{eqn:tot}). 
The theorem is an immediate consequence.
\end{proof}

\section{Algorithms for quadratic forms}
\label{sec:main}
In this section we describe the algorithm for computing the top
Betti numbers of a basic semi-algebraic set defined by quadratic forms.
We first describe an algorithm for computing index invariant
triangulations (see Definition \ref{def:iit}).

\begin{algorithm}[Index Invariant Triangulation]
\label{alg:triangulation}
\item[]
\item[{\sc Input}]
A set 
${\mathcal Q} = \{P_1,\ldots,P_\ell\} \subset \R[X_0,\ldots,X_k]$ 
where each $P_i$ is a quadratic form. We denote by
$P_{\mathcal Q} = (P_1,\ldots,P_\ell):\R^{k+1} \rightarrow \R^{\ell}$ the corresponding
quadratic map.
\item [{\sc Output}]
\item[]
\begin{enumerate}
\item
A semi-algebraic triangulation,
$$
h_{\mathcal Q}: \Delta_{\mathcal Q} \rightarrow \Omega_{\mathcal Q}
$$ 
of $\Omega_{\mathcal Q}$, 
which is compatible with the subsets $\Omega_{{\mathcal Q}'}$ for every
${\mathcal Q}' \subset {\mathcal Q}$, and
such that for any simplex $\sigma$ of $\Delta_{\mathcal Q}$, 
${\rm index}(\omega P)$, 
as well as the multiplicities of the eigenvalues of $\omega P$,
stay invariant as $\omega$ varies over $h_{\mathcal Q}(\sigma)$.
\item
For each simplex $\sigma$ of $\Delta_{\mathcal Q}$, 
a set of parametrized univariate representations,
$\{u_0(\sigma,\omega),\ldots,u_k(\sigma,\omega)\}$,
parametrized by $\omega \in h_{\mathcal Q}(\sigma)$,
such that the associated points in $\R^{k+1}$,
\[
\{e_0(\sigma,\omega),\ldots,e_k(\sigma,\omega)\},
\] 
form an orthornormal basis of $\R^{k+1}$, and 
\[
\{e_i(\sigma,\omega),\ldots,e_k(\sigma,\omega)\}
\] 
span the linear subspace of $\R^{k+1}$ on which the quadratic
form $\omega P_{\mathcal Q}$ is positive semi-definite (here $i  = {\rm index}(\omega P_{\mathcal Q})$).
\end{enumerate}

\item [{\sc Procedure}]
\item[]
\item[Step 1]
Let $\eps > 0$ be an infinitesimal and 
let $Z= (Z_1,\ldots,Z_\ell)$.
Also, let $M_{\mathcal Q}$ be 
the symmetric matrix corresponding to the quadratic form
(in $X_0,\ldots,X_k$) defined  by  
\[
\overline{P}_{\mathcal Q}(Z) = (1-\eps)(Z_{1}P_{1} +\cdots+ Z_{\ell}P_{\ell}) + \eps Q,
\]
where $Q = \sum_{i=0}^{k} i X_i^2$.
The entries
of $M_{\mathcal Q}$ depend linearly on $Z_1,\ldots,Z_\ell$, and
$\eps$. 
Compute the polynomial,
\[
F(Z,T) = \det(M_{\mathcal Q} + T\cdot I_{k+1})= 
T^{k+1} + C_{k}T^{k} + \cdots+ 
C_0,
\]
where each $C_i \in \R[\eps][Z_1,\ldots,Z_{\ell}]$ 
is a polynomial of degree at most $k+1$.
\item[Step 2]
Using Algorithm 11.1 in \cite{BPR03} (Elimination), 
compute a family of polynomials 
${\mathcal A}_{\mathcal Q} \subset \R[Z_1,\ldots,Z_\ell]$ such that 
such that 
for each $\rho \in {\rm Sign}({\mathcal A}_{\mathcal Q})$,
and $z \in \RR(\rho)$, 
the Thom encodings of the roots of $F(z,T)$ in $\R\la\eps\ra$, as well as the
the number of non-negative roots of $F(z,T)$ stay constant.

\item[Step 3]
Using the algorithm implicit in Theorem
 \ref{the:triangulation} (Triangulation),
compute a semi-algebraic triangulation,
$$
h_{\mathcal Q}: \Delta_{\mathcal Q} \rightarrow \Omega_{\mathcal Q},
$$
respecting the family  ${\mathcal A}_{\mathcal Q} \cup \cup_{i=1}^{\ell} 
\{Z_i\}$.

\item[Step 4]
For each simplex $\sigma$ of $\Delta_{\mathcal Q}$, 
let $\tau_0(\sigma),\ldots,\tau_k(\sigma)$ be the Thom encodings of 
the real roots,
$\lambda_0(\sigma,z) < \cdots < \lambda_k(\sigma,z)$  
of $F(z,T)$, for $z \in h_{\mathcal Q}(\sigma)$.
For $0 \leq i \leq k$, compute  using linear algebra a 
parametrized univariate representation 
\[
\bar{u}_i(Z,T) = (F,g_{i,0},\ldots,g_{i,k+1}),
\]
such that for
each $z \in h_{\mathcal Q}(\sigma)$ the point associated to
$\bar{u}_i(z,T)$ is the
eigenvector of $M_{\mathcal Q}(z)$ corresponding to the eigenvalue 
$\lambda_i(z)$.
Let $u_i = \lim_\eps \bar{u}_i$, and let $e_i(z) \in \R^{k+1}$ denote the
corresponding unit vector.
\end{algorithm}

\vspace{.1in}
\noindent
{\sc Complexity Analysis:}
The complexity of the algorithm is dominated by the complexity
of Step 3, which is $k^{2^{O(\ell)}}$.
\qedsymbol

\vspace{.1in}
\noindent
{\sc Proof of Correctness:}
\hide{
It follows from Descarte's rule of signs (see Remark 2.42, 
page 41 in \cite{BPR03}) that for any 
$z \in \Omega_{\mathcal Q}$,
${\rm index}(zP_{\mathcal Q})$ is equal to the number of sign variations in the
sequence $C_0(z),\ldots,C_{k}(z),+1$. Thus, the signs of the polynomials
${\mathcal A}_{\mathcal Q}= \{C_0,\ldots,C_k\}$ determine the index of $zP_{\mathcal Q}$. 
}
The quadratic form
$\overline{P}_{\mathcal Q}(z)$, and hence the matrix $M_{\mathcal Q}$, 
has $k+1$ distinct eigen values for each
$z\in \R^{\ell}$ and has the same index as 
$z P_{\mathcal Q}$. 
To prove the first part,  
replace $\eps$ in the definition of 
$\overline{P}_{\mathcal Q}(Z)$ and observe that the statement is true
$t=1$, since
the quadratic form $Q$ has distinct eigenvalues. Thus, the  set of $t$'s
for which   $\overline{P}_{\mathcal Q}(z)$ has $k+1$ distinct eigenvalues is non-empty,
constructible and contains a open subset, since the condition of having distinct
eigenvalues is a stable condition. Thus, there exists $\eps_0 > 0$, such that
for all $t\in (0,\eps_0)$, 
$\overline{P}_{\mathcal Q}(z)$ has $k+1$ distinct eigen-values, and hence
it is also the case for any infinitesimal $t$. 
The fact that $\overline{P}_{\mathcal Q}(z)$ has the same index as $zP_{\mathcal Q}$ 
follows from the fact that the  quadratic form $Q$ is positive semi-definite.
Thus, the eigenspaces corresponding to the eigenvalues of $M_{\mathcal Q}(z)$ are
all one-dimensional and thus the vectors $e_i'(\sigma,\omega)$ all well defined
(upto multiplication by $-1$).
Finally, note that since $e_0'(\sigma,\omega),\ldots,e_k'(\sigma,\omega)$ 
are orthonormal, so are $e_0(\sigma,\omega),\ldots,e_k(\sigma,\omega)$ for every
$\omega \in h_{\mathcal Q}(\sigma)$, and since 
$e_i'(\sigma,\omega),\ldots,e_k'(\sigma,\omega)$ span the non-negative eigenspace
of $\omega \overline{P}_{\omega}$, 
their images under the $\lim_\eps$ map will span the
non-negative eigenspace of  $\omega P_{\mathcal Q}$.
\qedsymbol

We now describe an algorithm for computing the
complexes $\check{\mathcal M}_{\mathcal Q}^{\bullet}$ 
described in the previous section.

\begin{algorithm}[Build Complex for Unions]
\label{alg:union}
\item[]
\item [{\sc Input}] 
\begin{enumerate}
\item An integer $\ell, 0 \leq \ell \leq k$.
\item
A quadratic map 
$P = (P_1,\ldots,P_s):\R^{k+1} \rightarrow \R^s$ given by $s$
homogeneous quadratic polynomials, $P_1,\ldots,P_s \in \R[X_0,\ldots,X_k].$
\end{enumerate}
\item [{\sc Output}]
\begin{enumerate}
\item
For each subset ${\mathcal Q} \subset {\mathcal P} =\{P_1,\ldots,P_s\}, 
\#{\mathcal Q} \leq {\ell + 2}$ 
a description of the
complex $\check{\mathcal M}_{\mathcal Q}$, consisting of a basis for each term of the
complex and matrices (in this basis) for the differentials.
\item
For each ${\mathcal Q}' \subset {\mathcal Q},$ with 
$\#{\mathcal Q} = \#{\mathcal Q}' + 1,$ matrices for the homomorphisms,
\[
\check{\phi}_{{\mathcal Q}',{\mathcal Q}}: \check{{\mathcal M}}^{\bullet}_{{\mathcal Q}'}
\rightarrow 
\check{{\mathcal M}}^{\bullet}_{\mathcal Q}.
\]
\end{enumerate}
\item [{\sc Procedure}]
\item[Step 1]
For each subset 
${\mathcal Q} = \{P_{i_1},\ldots,P_{i_{\ell + 2}}\} \subset {\mathcal P}$, with
$\#{\mathcal Q} = {\ell + 2}$,
let $P_{\mathcal Q}$ be the quadratic map corresponding to the subset 
${\mathcal Q}.$
Call Algorithm \ref{alg:triangulation}(Index Invariant Triangulation)
with input ${\mathcal Q}$.
\item[Step 2]
Construct the cell complex ${\mathcal C}(\Delta_{\mathcal Q})$ (following
its definition given in Section \ref{subsec:cellcomplex}).
\item[Step 3]
For each cell $C \in {\mathcal C}(\Delta_{\mathcal Q})$,
compute using 
the algorithm implicit in Theorem  \ref{the:triangulation} (Triangulation),
the cell complex ${\mathcal K}(C)$ and thus
obtain a description of ${\mathcal K}(B_{\mathcal Q})$.

\item[Step 4]
Compute the matrices corresponding to the differentials in the complex 
${\mathcal M}_{\mathcal Q}^{\bullet} = 
C^{\bullet}({\mathcal K}(B_{\mathcal Q})) $.
\item[Step 5]
For
${\mathcal Q}' \subset {\mathcal Q} \subset {\mathcal P}$ with
$\#{\mathcal Q} = \#{\mathcal Q}' + 1  < k$,
compute the matrices for the homomorphisms of complexes,
$$
\check{\phi}_{{\mathcal Q}',{\mathcal Q}}: \check{{\mathcal M}}^{\bullet}_{{\mathcal Q}'}
\rightarrow 
\check{{\mathcal M}}^{\bullet}_{\mathcal Q}.
$$
in the following way.

The simplicial complex ${\mathcal K}(B_{{\mathcal Q}'})$ 
is a subcomplex of 
${\mathcal K}(B_{\mathcal Q})$ by construction.
Compute the matrix for the restriction homomorphism,
$$
\phi_{{\mathcal Q},{\mathcal Q}'}:
C^{\bullet}({\mathcal K}(B_{{\mathcal Q}}))
\rightarrow 
C^{\bullet}({\mathcal K}(B_{{\mathcal Q}'})).
$$ 
and output the matrix for the dual homomorphism. 
\end{algorithm}

\vspace{.1in}
\noindent
{\sc Complexity Analysis:}
The complexity of Step 1 
$ 
\sum_{i=0}^{\ell+2} {s \choose i} k^{O(i)},
$
using the complexity of Algorithm \ref{alg:triangulation}.
The complexity of Step 2 is 
$ 
\sum_{i=0}^{\ell+2} {s \choose i} k^{2^{O(\min(\ell,s))}},
$
using the
complexity of the algorithm for triangulating semi-algebraic sets.
It follows from Proposition \ref{prop:complexity2}
that the complexities of all the remaining steps are
also bounded by 
$ 
\sum_{i=0}^{\ell+2} {s \choose i} k^{2^{O(\min(\ell,s))}}.
$
\qedsymbol

\vspace{.1in}
\noindent
{\sc Proof of Correctness:}
The correctness of the algorithm is a consequence of 
the correctness of Algorithm \ref{alg:triangulation} and 
Propositions
\ref{prop:iso2} and \ref{prop:commutative}.
\qedsymbol


Let $P_1,\ldots,P_s \in \R[X_0,\ldots,X_k]$ be homogeneous quadratic 
polynomials, and consider the set $S \subset \Sphere^k$ defined by,
$S = \{ x \in \Sphere^k \mid  P_1(x) \leq 0,\ldots,P_s \leq 0\}.$ 

We will also denote for $1 \leq i \leq s$, 
by $S_i$ the set defined by $\{x \in \Sphere^k \mid  P_i(x) \leq 0.\}$

Clearly, $S = \cap_{1 \leq i \leq s} S_i$.

\begin{algorithm} [Computing the highest $\ell$  Betti Numbers: 
the homogeneous case]
\label{alg:main}
\item[]
\item [{\sc Input}] A quadratic map 
$P = (P_1,\ldots,P_s):\R^{k+1} \rightarrow \R^s$ given by a set, 
${\mathcal P}=  \{P_1,\ldots,P_s \} \subset  \R[X_0,\ldots,X_k],$ of $s$
homogeneous quadratic polynomials.
\item [{\sc Output}] 
$b_{k}(S),\ldots, b_{k-\ell}(S)$, where $S$ is the set defined by
$$S = \bigcap_{P \in {\mathcal P}}
               \{x \in \Sphere^{k}\; \mid \; P(x) \leq 0 \}.
$$
\item [{\sc Procedure}]
\item[Step 1]
Using Algorithm \ref{alg:union} compute the 
truncated complex ${\mathcal D}_{\ell}^{\bullet,\bullet},$
i.e.
$$
\begin{array}{cccc}
{\mathcal D}_{\ell}^{p,q} &=& {\mathcal D}^{p,q}, &  \;\;0 \leq p \leq \ell+1,\;\; 
k-\ell-1 \leq q \leq k,  \\ 
                    &=& 0,             & \;\;\mbox{otherwise},
\end{array}
$$ 

\item[Step 2]
Compute using linear algebra, the ranks of 
\[
H^i({\rm Tot}^{\bullet}({\mathcal D}_\ell^{\bullet,\bullet})), 
\;\;k-\ell+1 \leq i \leq k. 
\]

\item[Step 3]
For each 
$i, \;\;k-\ell \leq i \leq k,$
output, 
$b_i(S) = {\rm rank}(H^i({\rm Tot}^{\bullet}({\mathcal D}_\ell^{\bullet,\bullet}))).$
\end{algorithm}

\vspace{.1in}
\noindent
{\sc Complexity Analysis:}
The number of algebraic operations is clearly bounded by
$ 
\sum_{i=0}^{\ell+2} {s \choose i} k^{2^{O(\min(\ell,s))}}
$
using the complexity analysis of 
Algorithm \ref{alg:union}.
\qedsymbol

\vspace{.1in}
\noindent
{\sc Proof of Correctness:}
The correctness of the algorithm is a consequence of the correctness of
Algorithm \ref{alg:union} and Theorem \ref{the:main}.
\qedsymbol

\begin{remark}
\label{rem:1}
Suppose that (using Notation from Algorithm \ref{alg:main})
${\mathcal P}' \subset {\mathcal P}$ and
$$
S' = \bigcap_{P \in {\mathcal P}'}
               \{x \in \Sphere^k\; \mid \; P(x) \leq 0 \},
$$
and letting ${\mathcal D}'^{\bullet,\bullet}_{\ell}$ denote the corresponding
complex for $S'$, it is clear from the definition that there is a 
homomorphism,
$\Phi_{{\mathcal P},{\mathcal P}'}: {\mathcal D}^{\bullet,\bullet}_\ell
\rightarrow {\mathcal D}'^{\bullet,\bullet}_\ell$ defined as follows.

For 
$$
\displaylines{\phi = \bigoplus_{{\mathcal Q} \subset {\mathcal P},\#{\mathcal Q} = p+1}
\phi_{\mathcal Q} \in  {\mathcal D}^{p,q}_\ell = 
\bigoplus_{{\mathcal Q} \subset {\mathcal P},\#{\mathcal Q} = p+1}
{\check{\mathcal M}}^{q}_{\mathcal Q}, \cr
\Phi_{{\mathcal P},{\mathcal P}'}(\phi) = 
\bigoplus_{{\mathcal Q} \subset {\mathcal P}',\#{\mathcal Q} = p+1}
\phi_{\mathcal Q}.
}
$$

Recall from \ref{eqn:fundamental} that there exists,
\[
\psi: {\mathcal D}_{\ell}^{\bullet,\bullet} 
\longrightarrow 
{\mathcal N}_{\ell}^{\bullet,\bullet}
\]
which induces an isomorphism,
$\psi^*:H^*(\Tot^{\bullet}({\mathcal D}_{\ell}^{\bullet,\bullet}) )
\longrightarrow 
H^*(\Tot^{\bullet}({\mathcal N}_{\ell}^{\bullet,\bullet})).$

Denoting by ${\mathcal N}'^{\bullet,\bullet}_\ell$ the 
(truncated) Mayer-Vietoris complex for
$S'$ and by 
$i_{{\mathcal P},{\mathcal P}'}: {\mathcal N}^{\bullet,\bullet}_\ell
\rightarrow
{\mathcal N}'^{\bullet,\bullet}_\ell
$ the inclusion homomorphism, we have the following commutative diagram.
\begin{diagram}
H^*(\Tot^{\bullet}({\mathcal D}^{\bullet,\bullet}_\ell))& 
\rTo^{\Phi_{{\mathcal P},{\mathcal P}'}^*}&
H^*(\Tot^{\bullet}({\mathcal D}'^{\bullet,\bullet}_\ell))\\
\dTo^{\psi^*}&&\dTo^{\psi'^*}\\
H^*(\Tot^{\bullet}({\mathcal N}^{\bullet,\bullet}_\ell))) &\rTo^{i_*} &
H^*(\Tot^{\bullet}({\mathcal N'}^{\bullet,\bullet}_\ell)))
\end{diagram}
Note that $H^*(\Tot^{\bullet}({\mathcal N}^{\bullet,\bullet})) \cong
H_*(S)$ and
$H^*(\Tot^{\bullet}({\mathcal N'}^{\bullet,\bullet}))\cong H_*(S').$

It is clear that Algorithm \ref{alg:main} can be easily modified to output the
complex ${\mathcal D}^{\bullet,\bullet}_\ell$, by outputting the matrices
corresponding to the vertical and horizontal homomorphisms in the
chosen bases.
Furthermore, given a subset 
${\mathcal P}' \subset {\mathcal P}$, Algorithm \ref{alg:main} 
can be made
to output both the complexes ${\mathcal D}^{\bullet,\bullet}_\ell$ and
${\mathcal D}'^{\bullet,\bullet}_\ell$ along with the matrices defining 
the homomorphism $\Phi_{{\mathcal P},{\mathcal P}'}$ with the same
complexity bounds.
\end{remark}

\section{The General Case}
\label{sec:general}
Let ${\mathcal P} = \{P_1,\ldots,P_s\} \subset \R[X_1,\ldots,X_k]$ with
${\rm deg}(P_i) \leq 2, 1 \leq i \leq s,$ and 
let $S \subset \R^k$ be the basic semi-algebraic set defined by
$P_1 \leq 0, \ldots, P_s \leq 0$. 

Let $\eps > 0$ be an infinitesimal,
and let $P_{s+1} = \eps\sum_{j=1}^{k} X_j^2 - 1.$ 
Let $\tilde{S} \subset \R\langle{\eps}\rangle^k$ 
be the basic semi-algebraic set defined by,
$P_1 \leq 0,\ldots, P_s \leq 0, P_{s+1} \leq 0.$
\begin{proposition}
\label{prop:perturb1}
The homology groups of $S$ and $\tilde{S}$ are isomorphic.
\end{proposition}
\begin{proof}
This is a consequence of the conical structure at infinity of 
semi-algebraic sets (see for instance Proposition 5.50 in \cite{BPR03})y. 
\end{proof}

Moreover, denoting by $P_i^h$ the homogenization of $P_i$,
and $\tilde{S}^h \subset \Sphere^k$ the
set defined by the system of quadratic inequalities,
\[
P^h_1 \leq 0,\ldots, P^h_s \leq 0, P^h_{s+1} \leq 0,
\]
on the unit sphere in $\R\langle {\eps}\rangle^{k+1}$
we have,

\begin{proposition}
\label{prop:perturb2}
For $0 \leq i \leq k,$
$b_i(\tilde{S}) = {1 \over 2}b_i(\tilde{S}^h).$
\end{proposition}

\begin{proof}
First observe that $\tilde{S}$ is bounded, and 
$\tilde{S}^h$ is the projection from the origin 
of the set $\tilde{S} \subset \{1\} \times \R\la\eps\ra^k$
onto the unit sphere in $\R\la\eps\ra^{k+1}$. Since, $\tilde{S}$ is bounded, 
the projection
does not intersect the equator and consists of two 
disconnected copies in the upper and lower
hemispheres, and each copy is homeomorphic to $\tilde{S}.$
\end{proof}

\begin{algorithm} 
[Computing the top $\ell$  Betti Numbers: the general case]
\label{alg:general}
\item[]
\item [{\sc Input}] A family of polynomials
$\{P_1,\ldots,P_s\}\subset  \R[X_1\ldots,X_k],$ with
${\rm deg}(P_i) \leq 2.$

\item [{\sc Output}] 
$b_{k-1}(S),\ldots, b_{k-\ell}(S)$, where $S$ is the set defined by
$$S = \bigcap_{P \in {\mathcal P}}
               \{x \in \R^{k}\; \mid \; P(x) \leq 0 \}.
$$
\item [{\sc Procedure}]
\item[Step 1]
Replace the family ${\mathcal P}$ by the family,
${\mathcal P}^h = \{ {P}^h_1,\ldots,{P}^h_s, P^h_{s+1} \}.$
\item[Step 2]
Using Algorithm \ref{alg:main} compute 
$b_k(\tilde{S}^h),\ldots,b_{k-\ell}(\tilde{S}^h).$

\item[Step 3]
Output $b_{k}(S) = {1\over 2}b_k(\tilde{S}^h), \ldots,
b_{k -\ell}(S) = {1\over 2}b_{k-\ell}(\tilde{S}^h).$
\end{algorithm}

\vspace{.1in}
\noindent
{\sc Proof of Correctness:}
The correctness of Algorithm \ref{alg:general} is a consequence of
Propositions \ref{prop:perturb1} and \ref{prop:perturb2} and the
correctness of Algorithm \ref{alg:main}.
\qedsymbol

\vspace{.1in}
\noindent
{\sc Complexity Analysis:}
The complexity of the algorithm is clearly 
\[
\sum_{i=0}^{\ell+2} {s \choose i} k^{2^{O(\min(\ell,s))}}
\]
from the complexity analysis of Algorithm \ref{alg:main}.
\qedsymbol

As in the case of Algorithm \ref{alg:main}(see Remark \ref{rem:1}) 
Algorithm \ref{alg:general} can be modified to 
output a complex whose associated total complex is
quasi-isomorphic to the truncated Mayer-Vietoris complex of $S$.
Furthermore, given a subset ${\mathcal P}' \subset {\mathcal P}$ defining
the set $S' \supset S$, Algorithm \ref{alg:general} can be made
to output the 
corresponding complexes of both sets, $S$ and $S'$,
as well as the homomorphism between them induced by inclusion.
Since this fact is important in certain applications we record it as
a separate theorem.

\begin{theorem}
\label{the:main2}
There exists an algorithm, which takes as input 
a family of polynomials
$\{P_1,\ldots,P_s\}\subset  \R[X_1\ldots,X_k],$ with
${\rm deg}(P_i) \leq 2,$
and a number $\ell \leq k$,
and outputs a complex ${\mathcal D}^{\bullet,\bullet}_\ell.$ 
The complex $\Tot^{\bullet}({\mathcal D}^{\bullet,\bullet}_\ell)$ is
quasi-isomorphic to $\C^\ell_{\bullet}(S)$, the truncated singular
chain complex of $S$,
where
$$S = \bigcap_{P \in {\mathcal P}}
               \{x \in \R^{k}\; \mid \; P(x) \leq 0 \}.
$$

Moreover, given a subset ${\mathcal P}' \subset {\mathcal P}$, 
with 
$$
S' = \bigcap_{P \in {\mathcal P}'}
               \{x \in \R^{k}\; \mid \; P(x) \leq 0 \}.
$$
the algorithm outputs both complexes ${\mathcal D}^{\bullet,\bullet}_\ell$ and
${\mathcal D}'^{\bullet,\bullet}_\ell$ (corresponding to the sets
$S$ and $S'$ respectively) along with the matrices defining 
a homomorphism $\Phi_{{\mathcal P},{\mathcal P}'},$ 
such that 
$\Phi_{{\mathcal P},{\mathcal P}'}^*: H_*(\Tot^\bullet({\mathcal D}^{\bullet,\bullet}_\ell)) \cong H_*(S) \rightarrow  H_*(S') \cong
H^*(\Tot^\bullet({\mathcal D'}^{\bullet,\bullet}_\ell))
$
is the homormorphism induced by the inclusion $i: S \hookrightarrow S'.$
The complexity of the algorithm is 
$ 
\sum_{i=0}^{\ell+2} {s \choose i} k^{2^{O(\min(\ell,s))}}.
$
\end{theorem}

\section{Hardness}
\label{sec:hardness}
In this section we show that the problem of 
computing the first few Betti numbers
of a semi-algebraic set defined by quadratic inequalities
is $\#$P-hard. Note that PSPACE-hardness of the problem of
counting the number of connected components for general semi-algebraic
sets were known before \cite{BC,Reif}
and the proofs of these
results extend easily to the  quadratic case.   

Recall that in the classical Turing machine model, a function taking its
values in ${\mathbb N}$ belongs to the class $\#$P, if there exists a 
polynomial time non-deterministic Turing machine, 
whose number of accepting paths on any
particular input is equal to the value of the function on that input 
\cite{Papadimitriou}.
A function $f$ is $\#$P-hard if every function in $\#$P can be reduced in
polynomial time to the computation of $f$. The problem of counting the number
of satisfying assignments to a Boolean formula is a $\#$P-hard
problem \cite{Papadimitriou}. 

The results of this section
complements the main algorithmic result of this paper 
proved in the last section, that the problem of computing the top few 
Betti numbers of such a set is in P.

\begin{theorem}
\label{the:hardness}
Given a family of polynomials ${\mathcal P} = \{P_1,\ldots,P_s\} \subset \R^k$, 
with $\deg(P) \leq 2, P \in {\mathcal P}$, as input,
the problem of computing $b_\ell(S)$, 
where $S \subset \R^k$ is defined
by $P_1 \geq 0,\ldots, P_s \geq 0$ and 
$\ell = O(\log k)$, is $\#$P-hard.
\end{theorem}

\begin{proof}
We first prove the case when $\ell = 0$ by a straightforward reduction from
the Boolean satisfiability problem. Clearly, if 
$\{C_1,\ldots,C_m\}$ is an instance of the Boolean satisfiability problem with
clauses $C_1,\ldots,C_m$ and variables $Z_1,\ldots,Z_n$, then the
number of satisfying assignments is equal to the number of connected
components of the set defined by,
$$
\displaylines{
X_1(X_1 - 1) \geq  0, X_1(1 - X_1) \geq  0,\ldots, 
X_n(X_n - 1) \geq  0, X_n(1 - X_n) \geq 0 \cr
\bar{C}_1 - 1 \geq 0,\ldots, \bar{C}_m -1  \geq 0,
}
$$
where $\bar{C}_i \in \R[X_1,\ldots,X_n]$ is the linear polynomial obtained
from the clause $C_i$ by syntactically substituting all the disjunctions by
additions, and for each $j, 1 \leq j \leq n,$
the literal $Z_j$ by the variable $X_j$, and the literal
$\neg Z_j$ by the expression $(1 - X_j)$.

For $\ell > 0$, we reduce to the case $\ell = 0$ using the following 
observation. 
Given a basic semi-algebraic set $S \subset \R^k$ defined by,
$P_1 \geq 0,\ldots, P_s \geq 0, \deg(P_i) \leq 2, 1 \leq i \leq s,$
one can define another basic semi-algebraic  set $S' \subset \R^{N}$
defined by $M$ polynomials inequalities with degrees $\leq 2$, 
with $M =s + 4m+1$ and $N = k+ 2m+1,$ 
where $m$ is the total number of monomials in the polynomials
$P_1,\ldots,P_s$. The number of monomials in the new system
is bounded by $10m$. 
The semi-algebraic set $S'$ has the homotopy type of the suspension $\Sigma S$.
Moreover, the description of $S'$ can be computed in polynomial time 
from  the description of $S$. It follows from the basic properties of the
suspensions (see \cite{Spanier})
that $b_1(\Sigma S) = b_0(S)$, which proves that computing
$b_1(S)$ is also $\#$P-hard. 
Iterating the construction, that is taking
suspensions of suspensions $\ell$ times, and noting that
$b_\ell(\Sigma^\ell(S)) = b_0(S),$ gives the result for 
$\ell = O(\log k)$ since the number of inequalities and variables
increases only polynomially in $s$ and $k$ for $\ell = O(\log k).$

We now describe the construction of the set $S'$. Introduce one new variable,
$X_0$ and consider the semi-algebraic set, $S'' \subset \R^{k+1}$ defined by,
$$
\displaylines{
(1 - X_0^2)P_1 \geq 0,\ldots, (1 - X_0^2)P_s \geq 0, \cr
1 - X_0^2 \geq 0.
}
$$
The set $S'' = ([-1,1] \times S) \cup H_1 \cup H_2$ where
$H_1$ and $H_2$ are the hyperplanes defined by $X_0 = 1$ and $X_0 = -1$
respectively. It is easy to see that $S''$ has the same homotopy type as the
suspension of $S$. However, the polynomials used in the description of $S''$ 
can have degrees as large as $4$. We show below that by introducing
new variables any quartic polynomial inequality can be reduced to 
a set of quadratic inequalities. The set $S'$ is defined as follows. 
Write each monomial $m$ of degree $> 2$ appearing in a polynomial 
used in the definition of $S''$, as a product of
two monomials, $m_1,m_2$,  each of degree at most $2$. Now replace
each ocurrence of $m$ in the inequalities used in the definition
of $S''$, by the quadratic monomial
$Y_{m,1}Y_{m,2}$, where $Y_{m,1},Y_{m,2}$ are new variables.
Also, add the inequalities
$$
\displaylines{
Y_{m,1} - m_1 \geq 0, m_1 - Y_{m,1} \geq 0,\cr 
Y_{m,2} - m_2 \geq 0, m_2 - Y_{m,2} \geq 0,
}
$$ 
to the set of inequalities.
Clearly, 
the number of monomials, the number of new variables and the number of
additional inequalities in the new system is bounded by 
$10m, 2m+1$ and $4m+1$ respectively, where
$m$ is the number of monomials in the original system. 
Finally, it is clear that the set $S'$ defined by the above inequalities
is a basic closed semi-algebraic set defined by polynomials
inequalities of degree at most $2$, which is homeomorphic to $S''$, and 
hence has the same homotopy type as $\Sigma S$. This completes the proof.
\end{proof}

\end{document}